\documentclass[letterpaper]{article} 

\usepackage[usenames,dvipsnames]{color}
\usepackage{bm}
\usepackage{graphicx,fancyhdr,lscape,dsfont,amsmath,amssymb} 
\usepackage{color, amsthm}
\usepackage{longtable}
\usepackage{arydshln}
\usepackage{psfrag}
\usepackage{multirow}
\usepackage{stmaryrd}
\usepackage[titletoc]{appendix}
\usepackage{graphicx}
\usepackage{ulem} 
\usepackage{mathrsfs}
\usepackage{bbm}
\usepackage{upgreek}
\usepackage{epstopdf}
\usepackage{booktabs}
\usepackage{tcolorbox}
\usepackage{appendix}
\usepackage{subcaption}

\allowdisplaybreaks

\newcommand{\vect}[1]{\vec{#1}} 
\newcommand{\tensor}[1]{  {\bm {#1}} } 

\newcommand{\e}{ {\vect{e}}}

\newcommand{\id}{ \mathds{1} }

\newcommand{\trace}[1]{ {\rm tr} \left[ \, {#1} \, \right] }

\newcommand{\deviatoric}[1]{ {\rm dev} \left[ \, {#1} \, \right] }
\newcommand{\divergence}[1]{ {\rm div} \left[ \, {#1} \, \right] }
\newcommand{\surfacedivergence}[2]{ {\rm div}_{#2} \left[ \, {#1} \, \right] }

\newcommand{\gradient}[1]{ {\rm \nabla} \left[ \, {#1} \, \right] }
\newcommand{\surfacegradient}[2]{ {\rm \nabla}_{#2} \left[ \, {#1} \, \right] }

\newcommand{\sym}[1]{ {\rm sym} \left[ \, {#1} \, \right] }

\newcommand{\bsigma}{\tensor{\sigma}}

\newcommand{\fpHad}{\mskip 3mu=\mskip-21mu \int}
\newcommand{\Cpv}{\mskip 3mu -\mskip-19mu \int}
\newcommand{\Guu}{ {\tensor{G}}_{uu} }
\newcommand{\Gup}{ {\tensor{G}}_{up} }
\newcommand{\Gpu}{ {\tensor{G}}_{pu} }
\newcommand{\Gpp}{ {\tensor{G}}_{pp} }

\newcommand{\Notin}{\mskip 3mu /\mskip-16mu \in}

\newcommand{\keywords}{{\em{Keywords }}}

\renewcommand{\em}[1]{\it{#1}}

\begin{document}

\title{A novel hydraulic fractures growth formulation.}

\author{  F. Fantoni$^{1,3}$, A. Salvadori$^{2,3}$
\\
\begin{small}
$^1$ DICATAM,  Universit\`a di Brescia, via Branze 43, 25123 Brescia, Italy
\end{small}
\\
\begin{small}
$^2$ DIMI,  Universit\`a di Brescia, via Branze 38, 25123 Brescia, Italy
\end{small}
\\
\begin{small}
$^3$ CeSiA,  Universit\`a di Brescia, via Branze 43, 25123 Brescia, Italy
\end{small}
 }

\date{}
\maketitle

\begin{abstract}
Propagation of a fluid-driven crack in an impermeable linear elastic medium under axis-symmetric conditions is investigated in the present work.
The fluid exerting the pressure  inside the crack is an incompressible Newtonian one and its front is allowed to lag behind the propagating fracture tip.
The tip cavity is considered as filled by fluid vapors under constant pressure having a negligible value with respect to the far field confining stress.
A novel algorithm is here presented, which is capable of tracking the evolution of both the fluid and the fracture fronts.
Particularly, the fracture tracking is grounded on a recent viscous regularization of the quasi-static crack propagation problem as a standard dissipative system. 
It allows a simple and effective approximation of the fracture front velocity by imposing Griffith's criterion at every propagation step.
Furthermore, for each new fracture configuration, a non linear system of integro-differential equations  has to be solved.
It arises from the non local elastic relationship existing between the crack opening and the fluid pressure, together with the non linear lubrication equation governing the flow of the fluid inside the fracture.    
\end{abstract}

\bigskip \noindent \keywords{Fracture mechanics, Hydraulic fracture, standard dissipative system.}

\bigskip \bigskip

\section{Introduction}
\label{Sec::Inroduction}

Hydraulic fractures identify a particular class of tensile fractures propagating in solids because of the pressure exerted by  injected viscous fluid under preexisting compressive stresses. They can be naturally present underground, as for magma transport in the Earth's crust \cite{Spence_Turcotte_1985}, or be man made fractures  which are created, for example, for stimulation of hydrocarbon reservoirs \cite{Economides_2000}, CO$_2$ sequestration \cite{Rudnicki_2000}, compensation grouting \cite{Mair_1995}, geothermal exploitation \cite{Legarth2005}, induced caving in mining \cite{Jeffrey_2000}.
Hydraulic fracture examples are very common in mechanobiology, too, since fluid-driven fractures are intimately connected to cells and tissues. They influence numerous cellular functions \cite{Casares2015,Lucantonio2017,Arroyo2017}, even from the very early stages of life \cite{Arroyo2019}.
Nowadays, the numerical simulation of fluid driven fractures continues to pose challenges and, depending upon the involved physical mechanisms, the complexity of the corresponding models may vary significantly \cite{Lecampion2018}.
The non local elastic relationship between the crack opening and the fluid pressure and the non linear lubrication equation governing the flow of the viscous fluid inside the cracks lead to a non linear system of inherently coupled integro-differential equations that must be solved for each new crack configuration. Literature shows that, even under the simplifying assumption of a plane strain or axial symmetric crack geometry, a rigorous mathematical solution of this moving-boundary problem is difficult to obtain \cite{DetournayReviewARFM2016}.
Furthermore, when the fluid front is allowed to lag behind the fracture tip, also the a priori unknown extent of such fluid lag has to be estimated.

The presence of a cavity behind the front of propagating cracks has been experimentally observed in low confining stress environments or in the presence of highly volatile fluids \cite{Medlin1984,Johnson1991,Groenenboom2001,Bunger2005}.
Precisely, because of the difficulty introduced by the presence of a fluid front lagging behind the leading fracture edge, most of previous and recent derived analytical solutions involving the hydromechanical coupling focus on zero-lag fractures (see, for example, \cite{Perkins1961,Geertsma1969,Abe1976,Spence1985,Huang1990,Carbonel999,Savitski2002,Garagash2005}).
Alternatively, different simplifying assumptions lessening the hydromechanical coupling have been speculated, such as the vanishing of pressure gradients inside the crack \cite{Khristianovic1955,Barenblatt1956,Jeffrey1989,Advani1997,Bui1996}.
Nevertheless, such hypothesis can be unrealistic since the fracturing fluid pressure field shows high gradients in the vicinity of the fracture front, thus leading to singular behavior, negative in sign, when the fluid and the crack edges coalesce in the LEFM context \cite{Spence1985,Lister1990,Desroches1994}.
In this sense, the lagging of the fluid front behind the advancing crack edge acquires a clear physical meaning, making the pressure field finite at the tip, without sustaining an infinite suction \cite{Garagash2006}.
Such a cavity is considered as filled by a fluid vapor having negligible value with respect to the far-field  confining stress if the surrounding medium can be modeled as impermeable, otherwise the pressure in the lag region is made equal to the pore pressure of the adjacent porous medium \cite{Rubin1993,Detournay2003}.
Fluid lags have been modeled in \cite{Lister1990,Garagash2000,Detournay2003} as extremely small compared to the fracture characteristic size, thus irrelevant for crack propagation.
A more general analysis of fluid lag evolution is accurately described in \cite{Garagash2006,Lecampion2007}   for plane strain fractures and in \cite{Bunger2007} for radial fractures at early stages of propagation. 
In \cite{Bunger2007} it was proved that  a characteristic time 
 scale exists, in the order of few seconds, in which cracks advance showing a fairly large lag behind the leading edge and the fluid pressure is significantly large compared to the confining stress. 
A further time scale, in the order of $10^5$ seconds, measures the time required by the fracture to evolve from a viscosity-dominated regime, when the main dissipative mechanism controlling crack growth is the viscous flow inside the fracture, to a toughness-dominated regime, when the major dissipation is due to the crack propagation and the consequent formation of new surfaces in the solid. 
 %
 %
 The size of the lag, therefore, declines with time and eventually vanishes.
 When the distance between the fluid and the crack fronts is large enough compared to the characteristic fracture length (early-time or modest confining stress), the stress field is dominated by the classical LEFM square-root asymptote (meaning the the fracture opening $w\sim \hat{x}^{1/2}$, where $\hat{x}$ measures the distance from the tip), because of the presence of a constant pressure in the lag zone. Hydromechanical coupling induces a $w\sim \hat{x}^{2/3}$ asymptotical behavior, which becomes dominant the smaller the lag  \cite{Adachi1994,Savitski2002,Lecampion2007}: this regime describes the large time response, or infinite confining stress solution, of the hydraulic fracture.

\bigskip
Several numerical algorithms have been designed to concurrently track two independent moving boundaries, following the evolution of the lag in time: the singular dislocation solutions,  either for the case of semi-infinite crack propagating at constant velocity \cite{Lister1990,Rubin1995,Garagash2000} or for the plane strain case \cite{Zhang2005}; finite elements \cite{Advani1997,Simoni2003,Schrefler2006}; XFEM \cite{Gordeliy2013,Gordeliy2013b,Gupta2014,Remij2015,Gupta2016,Gupta2018}.
{\em{The present work aims at establishing a set of governing equations that result in a novel numerical scheme, capable to accurately describe the evolution of a crack filled by a viscous Newtonian fluid in an infinite impermeable elastic medium. 
In the formulation, which stems from a standard dissipative picture of crack propagation in brittle materials \cite{SalvadoriJMPS08,SalvadoriJMPS10,SalvadoriCariniIJSS2011,SalvadoriFantoni2013,SalvadoriFantoniJECS2014,
SalvadoriFantoniIJSS2014,SalvadoriFantoniJMPS2016}, the fracturing fluid is allowed to lag behind the fracture front and the algorithm is capable of tracking both the moving crack and fluid edges concurrently.
While the fluid front velocity is computed in a closed form from the mass balance equation, the adopted crack tracking algorithm is grounded on a novel viscous regularization of the quasi-static crack propagation problem as a standard dissipative system recently presented in \cite{SalvadoriEtAlJMPS2019}. }}

\bigskip
A penny-shape propagating crack is used as a benchmark.
Despite the one-dimensional character of this axis-symmetric problem, the estimation of the velocity of the fluid-driven fracture at any given time is considered as challenging.
The proposed approach provides a simple and effective approximation of the fracture front velocity by imposing Griffith's criterion at every propagation step.
A hybrid scheme is employed: the elasto-hydrodynamic problem is solved implicitly, while the crack front is explicitly updated with the velocity computed at previous iteration.
%
%
To escape the limitation of propagation 
in pure mode I loading conditions, recent investigations on mixed mode hydraulic fracture propagation can be found in \cite{perkowska2017,wrobel2019}.

The paper is organized as follows:
Section \ref{Sec:HydromechanicalFormulation} is devoted to delineate the adopted hydro-mechanical formulation of the problem, how we model the response of the elastic domain and of the fluid.
Equations governing the two moving fronts are dealt with in the same Section, with particular emphasis to the viscous regularization technique adopted to approximate the fracture front velocity.
A numerical scheme for the propagation of a fluid-driven fracture in axis-symmetric conditions is detailed in Section \ref{Sec::Discretization}, validated in Section \ref{Sec::Benchmark} on a penny shape crack benchmark.
%

\section{ Hydro-mechanical problem formulation}
\label{Sec:HydromechanicalFormulation}

Consider an axis-symmetric fluid-driven fracture $\Gamma_w$ of radius $a(t)$ that, at initial time, has a given radius $a(0)>0$ and contains no fluid. Assume that, in its vicinity at least, the fracture is immersed in an impermeable, linear elastic medium with no limitations in stress and strain magnitude. In several applications, as underground $\rm CO_2$ storage, such a domain can be so wide compared to the fracture length to gain the notion of unboundedness: although we will take advantage of this assumption within this paper, this is not a mandatory condition.  

\begin{figure}[h!]
  \centering
  \includegraphics[width=7cm]{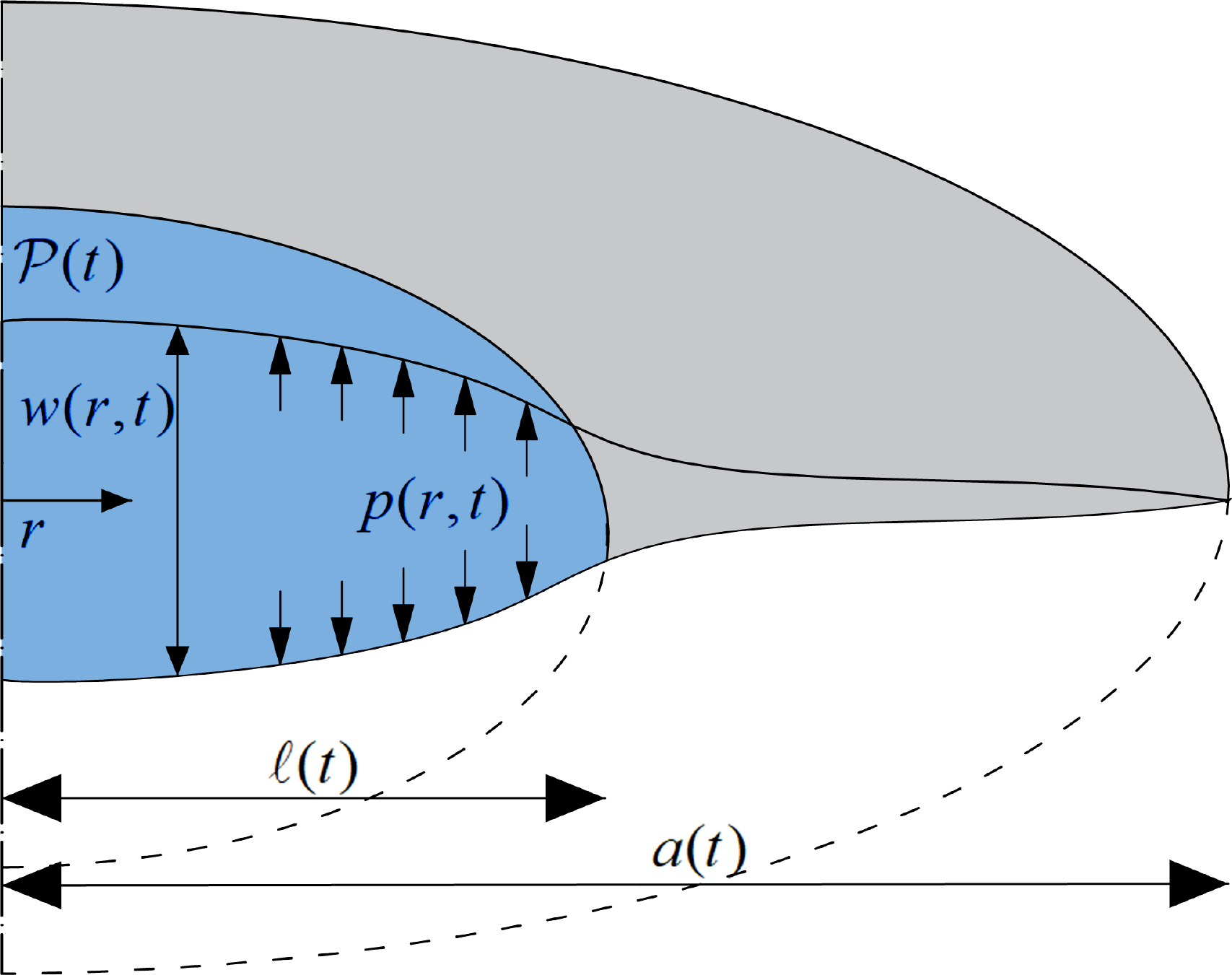}
  \caption{\it Penny shaped fracture of radius $a(t)$ filled by a Newtonian viscous fluid whose extent is $\ell(t)$.   }
  \label{Fig::fracture_draw}
\end{figure}

The driving force for the hydro-mechanical problem is assumed to be
a sufficiently smooth in time injection of an incompressible Newtonian fluid into the fracture. The injected fluid is assumed to occupy a region that is also axis-symmetric in time, with radius $\ell(t)\leq a(t)$, as depicted in Fig. \ref{Fig::fracture_draw}. 
The time-dependent fluid lag is defined as the difference between the crack and the fluid radii, namely $a(t)-\ell(t)$.
Denoting with $r$ the radial coordinate, for $r\leq\ell(t)$ the fracture walls are {\em{opened}} by a pressure equal to $p(r,t)$. 
In most cases of practical interest for hydraulic fracture, a superposition approach is taken and the problem is formulated in terms of a net pressure $p(r,t)$, the difference between a fluid pressure $p_f(r,t)$ and a far field compressive stress $\sigma_0$ orthogonal to the fracture plane.  
Since the surrounding medium is impermeable, the lag zone is filled by fluid vapor having a constant pressure that is assumed as negligible with respect to $\sigma_0$, thus $p(r,t)=-\sigma_0$ for $\ell(t)\leq r \leq a(t)$.

The unknown fields of the time dependent problem, formulated in the realm of small strains as usual for LEFM \cite{Rice1968}, are the crack opening $w(r,t)$, the net fluid pressure field $p(r,t)$, the fracture radius $a(t)$, and the fluid front position $\ell(t)$. They are all assumed to be smooth functions of time and space.

\subsection{Elastic response}
\label{subsec::UniboundedDomainElasticResponse}

\subsubsection{Differential form}

Denote with $\Gamma_u \subset \partial \Omega$ the Dirichlet part of the boundary, and with $\Gamma_p \subset \partial \Omega$ its complementary counterpart of Neumann type.  
Navier's mathematical formulation of the linear elastic fracture mechanics problem reads:
\begin{equation}
\label{Navier}
\begin{array}{lll}
 & 
 \divergence{ {\tensor{\sigma}} } + {\bar{\vect{f} }} = \vect{0}  & {\vect{x}} \in \Omega
 \cr
 & \vect{u} = \vect{\bar {u}} &   {\vect{x}} \in \Gamma_u 
 \cr
 & \tensor{\sigma} \, \vect{n} = \vect{\bar {p}} 
 & {\vect{x}} \in \Gamma_p 
 \cr
 & \tensor{\sigma} \, {\vect {n}}^+ = - p(r,t)  {\vect {n}}^+ & {\vect{x}} \in \Gamma_w^+ \cr
 & \tensor{\sigma} \, {\vect {n}}^- = - p(r,t)  {\vect {n}}^- & {\vect{x}}\in \Gamma_w^- \cr
\end{array}
\; ,
\end{equation}
where $\Gamma_w^+$ and $\Gamma_w^-$ are the upper and lower lips of the crack $\Gamma_w$.
In (isotropic) linear elastic fracture mechanics (LEFM), the Cauchy stress $\tensor{\sigma}$ is related to the small strain tensor ${\tensor{\varepsilon}}({\vect{u}})$ via the Helmholtz free energy $\Psi(\tensor{\varepsilon})$ as 
$$
\tensor{\sigma}({\tensor{\varepsilon}}) 
= 
\frac{\partial \Psi}{\partial \tensor{\varepsilon}}
= 
\lambda \, \trace{ \tensor{\varepsilon}}  \, \id 
+
2G \,  \tensor{\varepsilon}  
$$
where $\lambda$, G are the Lam{\`e}
constants and $\trace{\cdot}$ is the trace operator.
In order to obtain the weak form of (\ref{Navier}), we can proceed
as usual by multiplying for a test function ${\vect v}({\vect{x})}$ of a
suitable space V and integrating by parts in a
distributional sense \cite{QUARTERONI}.
The weak form of (\ref{Navier}) reads:

\begin{equation}
\label{WeakNavier}
find \;\; {{\vect{u}}} \in V : \; \; \; { \cal A} ({{\vect{u}}},{{\vect{v}}})
= {\cal F}({{\vect{v}}}) \; \; \;
\forall {{\vect{v}}} \in V
\end{equation}
where:
\begin{eqnarray*}
V &=& \left( H_{\Gamma_u}^1(\Omega) \right) ^ d 
\\
{ \cal A}
({{\vect{u}}},{{\vect{v}}}) 
&=&
\int_{\Omega} 
{\tensor{\varepsilon}}( {{\vect{v}}}) \; {:} \; 
{ \tensor{\sigma}({\tensor{\varepsilon}}( {\vect{u}}) ) }
\,
{\rm d} \Omega
\\
 {\cal F}({\vect{v}})
 &=&
 \int_{ \Omega} {{\vect{v}}}  \cdot {\bar {\vect{f}}} \, {\rm d} \Omega
 +
 \int_{\Gamma_p} {{\vect{v}}} \cdot {\bar {\vect{p}}} \, {\rm d} \Gamma
 +
 \int_{\Gamma_w} p(r,t) \; w(r,t)  \, {\rm d} \Gamma
 \;
 ,
\end{eqnarray*}
$d$ is the problem dimension and $\bar{{\vect{p}}} $, $\bar {{\vect{f}}} $
such that the integrals take sense, for instance $\bar {{\vect{f}}} \in
\left( L^2(\Omega) \right) ^ d $.
Existence and uniqueness of the weak solution of a problem are
usually stated by the Lax-Milgram theorem on a weak form of it
\cite{QUARTERONI}. 
It is easy to prove that ${\cal A}(\cdot,\cdot)$ is coercive,
i.e.
\begin{equation}
\label{coercivity}
 \exists \alpha>0 : { \cal A} ({{\vect{u}}}',{{\vect{u}}}') \geq
 \alpha \parallel {{\vect{u}}}' \parallel ^2  \; \; \forall {{\vect{u}}}'
 \in V
\end{equation}
Moreover, it is continuous,
i.e.
\begin{equation}
\label{continuity}
 \exists \gamma>0 : |{ \cal A} ({{\vect{u}}},{{\vect{v}}})| \leq
 \gamma \parallel {{\vect{u}}} \parallel \parallel {{\vect{v}}} \parallel  \; \;
 \forall {{\vect{u}}}, {{\vect{v}}} \in V
\end{equation}
From (\ref{coercivity}), (\ref{continuity}) and the continuity of
${\cal F}({{\vect{v}}})$ we get that the solution of the Navier's
differential problem \eqref{Navier} exists and it is unique provided that $ p(r,t) $ is given.

\subsubsection{Integral form}

The direct boundary integral formulation  \cite{RIZZO} of the elastic problem \eqref{Navier} rests on
Green's functions, gathered in matrices ${\tensor{G}}_{uu}$ and ${\tensor{G}}_{pu}$, which represent components ($i=1,2,3$) of the
displacement vector ${\vect{u}}$ at a point ${\vect{x}}$ and components
($i=1,2,3$) of the tractions vector ${\vect{p}}$ on a surface of normal
${{\vect{n}}}({{\vect{x}}})$ (i.e. ${\bsigma} ({{\vect{u}}}) {\vect{n}}$), due to a
unit force concentrated in space (at a point ${\vect{y}} \in \Omega$)
and acting on the unbounded elastic space $\Omega_\infty$ (that is
${\bar{\vect{f}}} = \delta({\vect{x}}-{\vect{y}}) {\vect{e}}_j$, $j=1,2,3$).
The expressions of the aforementioned Green's functions (also called
kernels) can be found in \cite{SalvadoriGray06,A-3DElAnInt}. 

Classical
mathematical passages (see e.g. \cite{LENBook}) lead to:
\begin{eqnarray*}
 {{\vect{u}}}({{\vect{x}}}) 
&+&
 \int_{\Gamma_p + \Gamma_u + \Gamma^+ + \Gamma^-}
 {\Gup} ({{\vect{x}}} - {{\vect{y}}}; {{\vect{l}}}({{\vect{y}}}) ) {{\vect{u}}}({{\vect{y}}}) \, {\rm d} {\Gamma_{y}}
=
\\
&&
 \int_{\Gamma_p + \Gamma_u + \Gamma^+ + \Gamma^-}
 {\Guu} ({{\vect{x}}} - {{\vect{y}}} ) {{\vect{p}}}({{\vect{y}}}) \, {\rm d} {\Gamma_{y}}
+
 \int_{\Omega}
 {\Guu} ({{\vect{x}}} - {{\vect{y}}} ) {\bar {{\vect{f}}}}({{\vect{y}}}) \, {\rm d} {\Gamma_{y}}
 \; ,
\qquad
{{\vect{x}}} \in \Omega
\end{eqnarray*}

\noindent
where ${{\vect{l}}}({{\vect{y}}})$ is the outward normal at ${{\vect{y}}} \in
\Gamma$. Since in small strains it holds \cite{ALBERTOTH}:
\begin{eqnarray}
\nonumber
 \int_{\Gamma^+ + \Gamma^-}
 {\Gup} ({{\vect{x}}} - {{\vect{y}}}; {{\vect{l}}}({{\vect{y}}}) ) {{\vect{u}}}({{\vect{y}}}) \, {\rm d} {\Gamma_{y}}
&=&
 \int_{\Gamma_w}
 {\Gup} ({{\vect{x}}} - {{\vect{y}}}; {{\vect{l}}}({{\vect{y}}}) ) {{\vect{w}}}({{\vect{y}}}) \, {\rm d} {\Gamma_{y}}
\\
\label{NoColl}
 \int_{\Gamma^+ + \Gamma^-}
 {\Guu} ({{\vect{x}}} - {{\vect{y}}} ) {{\vect{p}}}({{\vect{y}}}) \, {\rm d} {\Gamma_{y}} &=&
 {{\vect{0}}}
\end{eqnarray}
the so-called Somigliana's identity \cite{WATSON}
\begin{eqnarray}
\label{SomId}
\nonumber
{{\vect{u}}}({{\vect{x}}})
&+&
 \int_{\Gamma_p}
 {\Gup} ({{\vect{x}}} - {{\vect{y}}}; {{\vect{l}}}({{\vect{y}}}) ) {{\vect{u}}}({{\vect{y}}}) \, {\rm d} {\Gamma_{y}}
+
 \int_{\Gamma_u}
 {\Gup} ({{\vect{x}}} - {{\vect{y}}}; {{\vect{l}}}({{\vect{y}}}) ) {\bar {{\vect{u}}}}({{\vect{y}}}) \, {\rm d} {\Gamma_{y}}
\\ \nonumber
&+&
 \int_{\Gamma_w}
 {\Gup} ({{\vect{x}}} - {{\vect{y}}}; {{\vect{l}}}({{\vect{y}}}) ) {{\vect{w}}}({{\vect{y}}}) \, {\rm d} {\Gamma_{y}}
=
 \int_{\Gamma_u}
 {\Guu} ({{\vect{x}}} - {{\vect{y}}} ) {{\vect{p}}}({{\vect{y}}}) \, {\rm d} {\Gamma_{y}}
\\
&+&
 \int_{\Gamma_p}
 {\Guu} ({{\vect{x}}} - {{\vect{y}}} ) {\bar {{\vect{p}}}}({{\vect{y}}}) \, {\rm d} {\Gamma_{y}}
+
 \int_{\Omega}
 {\Guu} ({{\vect{x}}} - {{\vect{y}}} ) {\bar {{\vect{f}}}}({{\vect{y}}}) \, {\rm d} {\Gamma_{y}}
 \; ,
\qquad
{{\vect{x}}} \in \Omega
\end{eqnarray}
arises. Noteworthy Hong and Chen \cite{HONG-CHEN} remark that identity (\ref{NoColl}) represents
a mathematical degeneracy, whose consequence is the loss of
existence and uniqueness of the solution. It is clear that there
exists in the theory of elasticity a fundamental problem of a lack
of a general integral formulation for problems of an elastic body
with degenerate geometry that encloses no area or volume. 

To obtain an additional (dual) integral equation, Green's functions
(collected in matrices ${\tensor{G}}_{up}$ and ${\tensor{G}}_{pp}$) are
involved. They describe components ($i=1,2,3$) of the displacement
vector ${\vect{u}}$ at a point ${\vect{x}}$ and components ($i=1,2,3$) of the
tractions vector ${\vect{p}}$ on a surface of normal ${{\vect{n}}}({{\vect{x}}})$
(i.e. ${\bsigma} ({{\vect{u}}}) {{\vect{n}}}$), due to a unit relative
displacement concentrated in space (at a point ${\vect{y}} \in \Omega$)
and acting on the unbounded elastic space $\Omega_\infty$ (in
directions $j=1,2,3$). Applying the stress operator
on both sides of eq. (\ref{SomId}), with reference
to a normal ${{\vect{n}}}({{\vect{x}}})$ at a point ${{\vect{x}}} \in \Omega$, one
has:
\begin{eqnarray}
\label{TractEqOmega}
\nonumber
{{\vect{p}}}({{\vect{x}}}) 
&+&
 \int_{\Gamma_p}
 {\Gpp} ({{\vect{x}}} - {{\vect{y}}}; {{\vect{n}}}({{\vect{x}}}); {{\vect{l}}}({{\vect{y}}}) ) {{\vect{u}}}({{\vect{y}}})\, {\rm d} {\Gamma_{y}}
\\ \nonumber
&+&
 \int_{\Gamma_u}
 {\Gpp} ({{\vect{x}}} - {{\vect{y}}}; {{\vect{n}}}({{\vect{x}}}); {{\vect{l}}}({{\vect{y}}}) ) {\bar {{\vect{u}}}}({{\vect{y}}}) \, {\rm d} {\Gamma_{y}}
+
 \int_{\Gamma_w}
 {\Gpp} ({{\vect{x}}} - {{\vect{y}}}; {{\vect{n}}}({{\vect{x}}}); {{\vect{l}}}({{\vect{y}}}) ) {{\vect{w}}}({{\vect{y}}}) \, {\rm d} {\Gamma_{y}}
\\ \nonumber
&=&
 \int_{\Gamma_u}
 {\Gpu} ({{\vect{x}}} - {{\vect{y}}}; {{\vect{n}}}({{\vect{x}}}) ) {{\vect{p}}}({{\vect{y}}}) \, {\rm d} {\Gamma_{y}}
+
 \int_{\Gamma_p}
 {\Gpu} ({{\vect{x}}} - {{\vect{y}}}; {{\vect{n}}}({{\vect{x}}}) ) {\bar {{\vect{p}}}}({{\vect{y}}}) \, {\rm d} {\Gamma_{y}}
\\
&+&
 \int_{\Omega}
 {\Gpu} ({{\vect{x}}} - {{\vect{y}}}; {{\vect{n}}}({{\vect{x}}}) ) {\bar {{\vect{f}}}}({{\vect{y}}}) \, {\rm d} {\Gamma_{y}}
 \; ,
\qquad
{{\vect{x}}} \in \Omega
\end{eqnarray}
As for the adopted kernel notation, following \cite{POLIZZOTTO},
the first subscript of $\tensor{G}$ specifies the nature of the effect,
the second is associated with the quantity which is the
work-conjugate of the source causing that effect; in the argument
list, ${\vect{x}}$ and ${\vect{y}}$ denote field- and source-point,
respectively; the possible kernel dependence on the normal(s) is
explicitly indicated by a third (and a fourth) argument.

Let the classical bilinear form:
\begin{displaymath}
 {\cal A}({{\vect{u}}},{{\vect{p}}}) = \int_{\Gamma} {{\vect{u}}} \cdot {{\vect{p}}} \, {\rm d} {\Gamma_{x}}
\end{displaymath}
put into duality \cite{TONTI} the two spaces $U$ and $P$, with
${{\vect{u}}}
\in U$ and ${{\vect{p}}} \in P$. By applying Betti's theorem, the
following symmetry property (reciprocity relationships)
of the aforementioned Green's functions may be proved
\cite{HONG-CHEN}, \cite{SIRTORIMAIER}:
\begin{subequations}
\begin{eqnarray}
\label{Recipr1}
{\cal A}({{\vect{q}}}, {\Guu} {{\vect{p}}} ) &=& {\cal A}({{\vect{p}}}, {\Guu}
{{\vect{q}}} )
\\
\label{Recipr2}
{\cal A}({{\vect{u}}}, {\Gpp} {{\vect{v}}} ) &=& {\cal A}({{\vect{v}}}, {\Gpp}
{{\vect{u}}} )
\\
\label{Recipr3}
{\cal A}({{\vect{u}}}, {\Gpu} {{\vect{p}}} ) &=& {\cal A}({{\vect{p}}}, {\Gup}
{{\vect{u}}} )
\end{eqnarray}
\label{eq:Reciprocity}
\end{subequations}
Because of this property, equations
(\ref{SomId})-(\ref{TractEqOmega}) are called the dual equations
\cite{HONG-CHEN} for any point ${{\vect{x}}}
\in \Omega$.

All above introduced kernels are infinitely smooth functions
($C^\infty$) in their domain $D({{\vect{x}}}-{{\vect{y}}})$, which is the
whole space $\Re^d$ with exception of the origin. Therefore, all
integrals in equations (\ref{SomId})-(\ref{TractEqOmega}) are well
defined since ${{\vect{x}}} \Notin \Gamma$. If we take the point ${\vect{x}}$
to the boundary in a limit process, singularities will arise into
integrals in equations (\ref{SomId})-(\ref{TractEqOmega}). Those
singularities have been extensively investigated during the last
years (see e.g. \cite{BREBBIA,DILIMONEGATO,GuiggianiBook}): their singularity-orders are collected in table
(\ref{Table:Singul}).

\begin{table}[ht]
{\centering
 \begin{tabular}{| c | c | c | c | } \hline
     kernel   & \multicolumn{2}{|c|}{Asymptotical behavior}  & Denomination \\
              & \multicolumn{2}{|c|}{ when ${{\vect{r}}} \rightarrow {{\vect{0}}}$}  &  \\  \hline
              &   2D  &  3D & \\ \hline
      $\Guu$  &  $O( \log(r))$ &  $O( r^{-1})$ &  Weak singularity (integrable) \\ \hline
      $\Gup$,$\Gpu$  &  $O(r^{-1}) $ &  $O(r^{-2})$ &  Strong singularity \\ \hline
      $\Gpp$  &  $O( r^{-2}) $ &  $O( r^{-3})$ &  Hyper singularity \\ \hline
 \end{tabular}
\par}
\caption{ \it Kernels and their singularities. }
\label{Table:Singul}
\end{table}

Strongly singular kernels $\Gup$ and $\Gpu$ generate free terms
($c({{\vect{x}}})$ and $d({{\vect{x}}})$ respectively) in the limit process
\cite{Gray2, Guiggiani}; free terms are equal to 1/2 for
smooth boundaries. ``Integrals" involving strongly singular kernels
must be understood in their distributional nature of Cauchy
Principal Value \cite{RjasanowSteinbachBook}. Similarly, the hyper singular
kernel $\Gpp$ must be understood in its distributional nature of
Hadamard's finite part \cite{DILIMONEGATO}. The following
boundary integral equations eventually arise:
\begin{subequations}
\begin{eqnarray}
\nonumber
c({{\vect{x}}}) {{\vect{u}}}({{\vect{x}}}) 
&+&
 \Cpv_{\Gamma_p}
 {\Gup} ({{\vect{x}}} - {{\vect{y}}}; {{\vect{l}}}({{\vect{y}}}) ) {{\vect{u}}}({{\vect{y}}}) {\rm d} \, {\rm d} {\Gamma_{y}}
+
 \Cpv_{\Gamma_u}
 {\Gup} ({{\vect{x}}} - {{\vect{y}}}; {{\vect{l}}}({{\vect{y}}}) ) {\bar {{\vect{u}}}}({{\vect{y}}}) \, {\rm d} {\Gamma_{y}}
\\ \nonumber
&+&
 \Cpv_{\Gamma_w}
 {\Gup} ({{\vect{x}}} - {{\vect{y}}}; {{\vect{l}}}({{\vect{y}}}) ) {{\vect{w}}}({{\vect{y}}}) \, {\rm d} {\Gamma_{y}}
=
 \int_{\Gamma_u}
 {\Guu} ({{\vect{x}}} - {{\vect{y}}} ) {{\vect{p}}}({{\vect{y}}}) \, {\rm d} {\Gamma_{y}}
\\
&+&
 \int_{\Gamma_p}
 {\Guu} ({{\vect{x}}} - {{\vect{y}}} ) {\bar {{\vect{p}}}}({{\vect{y}}}) \, {\rm d} {\Gamma_{y}}
+
 \int_{\Omega}
 {\Guu} ({{\vect{x}}} - {{\vect{y}}} ) {\bar {{\vect{f}}}}({{\vect{y}}}) \, {\rm d} {\Gamma_{y}}
 \; ,
\qquad \qquad
{{\vect{x}}} \in \Gamma
\label{DisplEq}
\\
\nonumber
d({{\vect{x}}}) {{\vect{p}}}({{\vect{x}}}) 
&+&
 \fpHad_{\Gamma_p}
 {\Gpp} ({{\vect{x}}} - {{\vect{y}}}; {{\vect{n}}}({{\vect{x}}}); {{\vect{l}}}({{\vect{y}}}) ) {{\vect{u}}}({{\vect{y}}}) \, {\rm d} {\Gamma_{y}}
\\ \nonumber
&+&
 \fpHad_{\Gamma_u}
 {\Gpp} ({{\vect{x}}} - {{\vect{y}}}; {{\vect{n}}}({{\vect{x}}}); {{\vect{l}}}({{\vect{y}}}) ) {\bar {{\vect{u}}}}({{\vect{y}}}) \, {\rm d} {\Gamma_{y}}
+
 \fpHad_{\Gamma_w}
 {\Gpp} ({{\vect{x}}} - {{\vect{y}}}; {{\vect{n}}}({{\vect{x}}}); {{\vect{l}}}({{\vect{y}}}) ) {{\vect{w}}}({{\vect{y}}}) \, {\rm d} {\Gamma_{y}}
\\ \nonumber
&=&
 \Cpv_{\Gamma_u}
 {\Gpu} ({{\vect{x}}} - {{\vect{y}}}; {{\vect{n}}}({{\vect{x}}}) ) {{\vect{p}}}({{\vect{y}}}) \, {\rm d} {\Gamma_{y}}
+
 \Cpv_{\Gamma_p}
 {\Gpu} ({{\vect{x}}} - {{\vect{y}}}; {{\vect{n}}}({{\vect{x}}}) ) {\bar {{\vect{p}}}}({{\vect{y}}}) \, {\rm d} {\Gamma_{y}}
\\
&+&
 \Cpv_{\Omega}
 {\Gpu} ({{\vect{x}}} - {{\vect{y}}}; {{\vect{n}}}({{\vect{x}}}) ) {\bar {{\vect{f}}}}({{\vect{y}}}) \, {\rm d} {\Gamma_{y}}
 \; ,
\qquad \qquad
{{\vect{x}}} \in \Gamma
\label{TractEq}
\end{eqnarray}
\label{eq:IntegralEqs}
\end{subequations}

Equation (\ref{DisplEq}), often referred to as displacement
equation, is the starting point for the numerical approximation of
elasticity problems via the boundary element method (BEM) by the
collocation technique \cite{BREBBIA,WATSON}. In
modeling crack problems, identity (\ref{NoColl}) represents an
insurmountable mathematical difficulty in applying the boundary
element method via the collocation technique making use of the
displacement equation only (see e.g. \cite{CRUSE, ALROOK}).
Some special techniques have been devised to overcome this
mathematical degeneracy. Among others, the special Green's
functions methods \cite{CRUSE2}, the zone method \cite{LIGGET-LIU}
and the Dual BEM \cite{PORTELA-ALIABADI}. This last technique makes
use of (\ref{TractEq}), often referred to as traction equation.

As an alternative approach, a Galerkin scheme for the BIE (Boundary Integral Equation)
formulation has been proposed (see e.g. \cite{LENBook,NEDELEC-1,SirtoriMeccanica1979,HARTMANN}). To obtain such a formulation
eq. (\ref{DisplEq}) is imposed on the Dirichlet boundary
$\Gamma_u$, whereas eq. (\ref{TractEq}) is imposed both on the
Neumann boundary $\Gamma_p$ and on the two boundaries $\Gamma_w^+$,
$\Gamma_w^-$. In an operatorial form:
\begin{eqnarray}
\label{Overview:IntEq}
\left[ \begin{array}{ccc}
\int_{\Gamma_{u}}\Guu  [.] \, {\rm d} {\Gamma_{y}}   &
- \Cpv_{\Gamma_p} \Gup  [.] \, {\rm d} {\Gamma_{y}}   &
- \Cpv_{\Gamma_p} \Gup  [.] \, {\rm d} {\Gamma_{y}}    \\
- \Cpv_{\Gamma_u} \Gpu  [.] \, {\rm d} {\Gamma_{y}}   &
  \fpHad_{\Gamma_p} \Gpp  [.] \, {\rm d} {\Gamma_{y}}    &
  \fpHad_{\Gamma_p} \Gpp  [.] \, {\rm d} {\Gamma_{y}}    \\
- \Cpv_{\Gamma_u} \Gpu  [.] \, {\rm d} {\Gamma_{y}}   &
  \fpHad_{\Gamma_p} \Gpp [.] \, {\rm d} {\Gamma_{y}} &
  \fpHad_{\Gamma_p} \Gpp  [.] \, {\rm d} {\Gamma_{y}}
\end{array}
\right]
              \left[ \begin{array}{c}
               {{\vect{p}}}   \\  {{\vect{u}}}  \\  {{\vect{w}}}
                \end{array}
\right]
%
+
\left[ \begin{array}{c}
                 {{\vect{0}}}    \\     {{\vect{0}}}   \\
                 {{\vect{p}}}
                 \end{array}
\right]
=
\left[ \begin{array}{c}
                 {{\vect{f}}}^u    \\     {{\vect{f}}}^p   \\
                 {{\vect{f}}}^w
                 \end{array}
\right] \begin{array}{c}
              \mbox{on $\Gamma_u$} \\  \mbox{on $\Gamma_p$} \\  \mbox{on $\Gamma_w$} \\
              \end{array}
\end{eqnarray}
Vectors $ {{\vect{f}}}^i \; , \; i=u,p,w$, that gather all data (i.e. $
{\bar{
{\vect{p}}}}$, $ {\bar{{\vect{u}}}}$, $ {\bar{{\vect{f}}}}$), are as follows:
\begin{eqnarray*}
{ {{\vect{f}}}}^u \left( {{\vect{x}}} \right) &:=& c({{\vect{x}}}){
\bar {{\vect{u}}}} -
\int_{\Gamma_p} \Guu   { \bar{{\vect{p}}} } \, {\rm d} {\Gamma_{y}} +
\Cpv_{\Gamma_u} \Gup   { \bar{{\vect{u}}} } \, {\rm d} {\Gamma_{y}} -
  \int_{\Omega} \Guu   { \bar{{\vect{f}}} } \, {\rm d} {\Gamma_{y}}
\\
{ {{\vect{f}}}}^p \left( {{\vect{x}}} \right) &:=& - d({{\vect{x}}})\bar{{\vect{p}}} +
  \Cpv_{\Gamma_p} \Gpu   {\bar{{\vect{p}}}} \, {\rm d} {\Gamma_{y}} -
\fpHad_{\Gamma_u} \Gpp   {\bar{{\vect{u}}}} \, {\rm d} {\Gamma_{y}} +
    \Cpv_{\Omega} \Gpu   {\bar{{\vect{f}}}} \, {\rm d} {\Gamma_{y}}
\\
 {{\vect{f}}}^{d} \left( {{\vect{x}}} \right)  &:=&
 \Cpv_{\Gamma_p} \Gpu  {\bar{{\vect{p}}}} \, {\rm d} {\Gamma_{y}} -
 \fpHad_{\Gamma_u} \Gpp  {\bar{{\vect{u}}}} \, {\rm d} {\Gamma_{y}} +
    \Cpv_{\Omega} \Gpu  {\bar{{\vect{f}}}} \, {\rm d} {\Gamma_{y}}
\end{eqnarray*}
%
%
The problem \eqref{Overview:IntEq} admits a variational representation: displacements
${{\vect{u}}}$ on the free boundary $\Gamma_p$, tractions ${{\vect{p}}}$ on
the constrained boundary $\Gamma_u$ and the crack opening
displacements ${{\vect{w}}}$ on the boundary $\Gamma_w$ which solve the
problem are characterized by the stationarity of the following
quadratic functional:
\begin{eqnarray}
\label{Overview:eq:Func}
 \Psi [ {{\vect{u}}}, {{\vect{p}}}, {{\vect{w}}} ]
 &=& 
  \frac{1}{2}  \left\{ \;  \int_{\Gamma_{u}}  {{\vect{p}}} ({{\vect{x}}})  \cdot \int_{\Gamma_{u}}\Guu  \;  {{\vect{p}}} ({{\vect{y}}}) \, {\rm d} {\Gamma_{y}} \, {\rm d} {\Gamma_{x}}
               \right.
 \\
 &-& \nonumber
 \int_{\Gamma_{u}}  {{\vect{p}}} ({{\vect{x}}}) \cdot \Cpv_{\Gamma_{p}}\Gup   \; {{\vect{u}}} ({{\vect{y}}}) \, {\rm d} {\Gamma_{y}} \, {\rm d} {\Gamma_{x}}
 -
 \int_{\Gamma_{u}} {{\vect{p}}} ({{\vect{x}}})  \cdot \Cpv_{\Gamma_w}\Gup  \; {{\vect{w}}} ({{\vect{y}}}) \, {\rm d} {\Gamma_{y}} \, {\rm d} {\Gamma_{x}}
 \\
 &-& \nonumber
 \int_{\Gamma_{p}}  {{\vect{u}}} ({{\vect{x}}})  \cdot \Cpv_{\Gamma_{u}}\Gpu   \;   {{\vect{p}}} ({{\vect{y}}}) \, {\rm d} {\Gamma_{y}} \, {\rm d} {\Gamma_{x}}
 +
 \int_{\Gamma_{p}} {{\vect{u}}} ({{\vect{x}}})  \cdot \fpHad_{\Gamma_{p}}\Gpp  \;  {{\vect{u}}} ({{\vect{y}}}) \, {\rm d} {\Gamma_{y}} \, {\rm d} {\Gamma_{x}}
 \\
 &+& \nonumber
 \int_{\Gamma_{p}}  {{\vect{u}}} ({{\vect{x}}})  \cdot \fpHad_{\Gamma_w}\Gpp  \;    {{\vect{w}}} ({{\vect{y}}}) \, {\rm d} {\Gamma_{y}} \, {\rm d} {\Gamma_{x}}
 -
 \int_{\Gamma_w} {{\vect{w}}} ({{\vect{x}}})  \cdot \Cpv_{\Gamma_{u}}\Gpu  \;  {{\vect{p}}} ({{\vect{y}}}) \, {\rm d} {\Gamma_{y}} \, {\rm d} {\Gamma_{x}}
 \\
 &+& \nonumber
\left.
 \int_{\Gamma_w}  {{\vect{w}}} ({{\vect{x}}})  \cdot \fpHad_{\Gamma_{p}}\Gpp  \;   {{\vect{u}}} ({{\vect{y}}}) \, {\rm d} {\Gamma_{y}} \, {\rm d} {\Gamma_{x}}
+
 \int_{\Gamma_w}  {{\vect{w}}} ({{\vect{x}}})  \cdot \fpHad_{\Gamma_{w}}\Gpp  \;    {{\vect{w}}} ({{\vect{y}}}) \, {\rm d} {\Gamma_{y}} \, {\rm d} {\Gamma_{x}}
\right\}
\\
&-& \nonumber
\int_{\Gamma_{u}}  {{\vect{p}}} ({{\vect{x}}})  \cdot \;   {{{\vect{f}}}}^{u} ({{\vect{x}}}) \, {\rm d} {\Gamma_{x}}
-
\int_{\Gamma_{p}}  {{\vect{u}}} ({{\vect{x}}})  \cdot \;   {{{\vect{f}}}}^{p} ({{\vect{x}}}) \, {\rm d} {\Gamma_{x}}
-
\int_{\Gamma_w}  {{\vect{w}}} ({{\vect{x}}})  \cdot \;  
   \left( 
   {{{\vect{f}}}}^{w} ({{\vect{x}}}) - {{{\vect{p}}}} ({{\vect{x}}})
   \right)
   \, {\rm d} {\Gamma_{x}}
\end{eqnarray}

\noindent
The stationary point is a {\em saddle-point} (minimum with respect
to $ {{\vect{p}}}$ and maximum with respect to both $ {{\vect{u}}}$ and $
{{\vect{w}}}$ ). The functional (\ref{Overview:eq:Func}) is the
counterpart of the Hu-Washizu functional for the boundary integral
formulation \cite{POLIZZOTTO}.
Existence and uniqueness of the weak solution are guaranteed by the
regularity properties of the involved integral operators in
suitable Hilbert spaces (see e.g. \cite{ALETESI},
\cite{WENDLAND}).

\subsubsection{Integral equations for a plane crack in an infinite medium}
\label{sec:plcr}

As a special case of the general framework depicted so far,  
consider a pressurized crack lying in the plane $\e_2 \times
\e_3$, embedded in an infinite elastic medium with negligible
volume forces. By the considered hypotheses, 
$ \Gamma_u = \Gamma_p= \emptyset$ and 
$\bar{{\vect{f}}} = {{\vect{0}}}$. Consider further ${\vect{n}}^+ = -\e_1$ and ${{\vect{p}}}^+ = p({\vect{x}},t) \, \e_1$ .
Problem (\ref{Overview:IntEq}) becomes therefore:
\begin{equation}
\label{CrackEq}
 \fpHad_{\Gamma_w}
 {\Gpp} ({{\vect{x}}} - {{\vect{y}}}; -\e_1; -\e_1 ) \; {{\vect{w}}}({{\vect{y}}}) \, {\rm d} {\Gamma_{y}}
 = - p({\vect{x}},t) \, \e_1
 \; ,
\qquad \qquad {{\vect{x}}} \in \Gamma_w
\end{equation}
where kernel $\Gpp$ yields:
\begin{eqnarray*}
 \Gpp \left( {{\vect{r}}} ; -\e_1; -\e_1 \right)
 =
  \frac{G \nu}{4 \pi (1-\nu)} \frac{1}{r^3}
  \left\{
   2 \, (\e_1 \otimes \e_1)
+
 3
 \frac{{{\vect{r}}} \otimes {{\vect{r}}}}{r^2}  \;
 +
 \frac{(1-2\nu)}{\nu}
  {\id}
  \right\}
\end{eqnarray*}
If one writes equation (\ref{CrackEq}) by components (${{\vect{x}}} \in
\Gamma_w$):
\begin{equation}
\label{CrackEqByComp}
 \frac{G \nu}{4 \pi (1-\nu)}  \fpHad_{\Gamma_w}
 \frac{1}{r^3}
 \left[
 \begin{array}{c}
 \frac{w_1({{\vect{y}}})}{ \nu} \cr
 \left( \frac{1-2\nu}{\nu} + 3 \frac{r_2^2}{r^2} \right) w_2({{\vect{y}}}) + 3 \frac{r_2 \, r_3}{r^2} w_3({{\vect{y}}}) \cr
 3 \frac{r_2 \, r_3}{r^2} w_2({{\vect{y}}}) + \left( \frac{1-2\nu}{\nu} + 3 \frac{r_3^2}{r^2} \right) w_3({{\vect{y}}})
 \end{array}
 \right]
 \;
 {\rm\,d} \Gamma_{y}
 = \left[
 \begin{array}{c}
 - p({\vect{x}},t) \cr 0 \cr 0
 \end{array}
 \right]
\end{equation}
one notes that a solution for (\ref{CrackEq}) is:
\begin{equation}
 \label{CrackEqSol}
 w_2({{\vect{y}}})=w_3({{\vect{y}}})=0 \, , \qquad w_1({{\vect{y}}}) \quad \mbox{s. t.} \quad
 \fpHad_{\Gamma_w} \frac{w_1({{\vect{y}}})}{r^3}  \; {\rm\,d} \Gamma_{y}
 = - 8 \pi \, \frac{1-\nu^2}{E}  \, p({\vect{x}},t)
\end{equation}
%
It is
worthy here to stress that even if
\begin{displaymath}
 \frac{w_1({{\vect{y}}})}{r^3} \ge 0 \qquad \mbox {when } {{\vect{x}}}, {{\vect{y}}} \in \Gamma_w
\end{displaymath}
its finite part of Hadamard can be negative, as it is expected in
equation (\ref{CrackEqSol}-b). For a larger comprehension, see
\cite{SalvadoriIJNME07,SalvadoriIJNME10}.

\subsubsection{Weight functions and unbounded domains}
\label{sec:wfud}

For the penny shaped fracture of radius $a(t)$ filled by a Newtonian viscous fluid depicted in Fig. \ref{Fig::fracture_draw},
acting in an
unbounded domain, in view of the simple geometry of the crack the elastic response can be expressed through appropriate weight functions, which provide the crack opening $w(r,t)$ in terms of the net pressure $p(r,t)$ \cite{Sneddon1951} as
\begin{equation}
\label{eq:wSneddon}
w(r,t)=\frac{8 a}{\pi E'}\int_{r/a}^1 \frac{\xi}{\sqrt{\xi^2-(r/a)^2}}
\int_{0}^{1}\frac{x\,p(x\xi a,t)}{\sqrt{1-x^2}}\, {\rm d}x\, {\rm d}\xi,
\end{equation}
where  $E'=\frac{E}{1-\nu^2}$ is the plane strain modulus expressed in terms of Young modulus $E$ and Poisson ratio $\nu$ of the host material. Note that eq. \eqref{eq:wSneddon} is less general than \eqref{CrackEqSol}.
In order to express equation (\ref{eq:wSneddon}) in a more useful formalism, variable substitutions can be performed eventually leading from \eqref{eq:wSneddon} to
%
%
\begin{equation}
\label{eq:wNostra}
w(r,t)=\frac{8 }{\pi E'}\int_{r}^a \frac{1}{\sqrt{z^2-r^2}}
\int_{0}^{z}\frac{y\,p(y,t)}{\sqrt{z^2-y^2}}\, {\rm d}y \, {\rm d}z
\; .
\end{equation}
%

%
\subsection{Fluid response}
\label{subsec::LubricationEquation}

\subsubsection{Mass balance}

The {\em{mass balance}} for a convecting volume $\mathcal{P}(t)$ reads \cite{GurtinFriedAnand}
\begin{equation}
\frac{d}{dt}\int_{\mathcal{P}(t)}\rho\, {\rm d}v
=
- 
\int_{\partial\mathcal{P}(t)}
\vect{h}\cdot\vect{n}
\, {\rm d}a
+
\int_{\mathcal{P}(t)} s \, {\rm d} v 
\; ,
\label{eq:MassBalance}
\end{equation}
where $\rho$ is the density (mass per unit volume), $\vect{h}$ is the flow of mass across the convecting boundary $\partial\mathcal{P}(t)$, $\vect{n}$ is the outward unit normal to $\partial\mathcal{P}(t)$.
%
It is here assumed that the driving force, namely the injection of fracturing fluid, is accounted for by means of a sufficiently smooth in time and in space source term $s=s({\vect{x}},t)$.
%
Exploiting Reynold's transport theorem, equation (\ref{eq:MassBalance}) can be written in a localized form as
\begin{equation}
\frac{\partial \rho}{\partial t} +  \divergence{\rho\vect{v}}=s
\; ,
\label{eq:MassBalancaLocalizedForm}
\end{equation}
where $\vect{v}$ is the fluid velocity and symbol $\divergence{\cdot}$ stands for the divergence operator.
Considering an incompressible and homogeneous fluid, for which ${\partial \rho}/{\partial t}=0$ and $\gradient{\rho}=\vect{0}$, mass balance (\ref{eq:MassBalancaLocalizedForm}) eventually simplifies as
 \begin{equation}
  \divergence{\vect{v}} = s/\rho
\label{eq:MassBalancaLocalizedForm2}
\end{equation}
in the current configuration, $ ({\vect{x}},t) \in \mathcal{P}(t) \times [0,T] $.

\subsubsection{Momentum balance}

Assuming negligible inertia and body forces, the
balance of momentum in the current configuration, i.e. $ ({\vect{x}},t) \in \mathcal{P}(t) \times [0,T] $, reads
\begin{equation}
\divergence{\tensor{\sigma}}=\vect{0}
\; .
\label{eq:MomentumBalance}
\end{equation}
The Cauchy stress tensor is further decomposed in its spherical and deviatoric part, with pressure $ p = -\trace{\tensor{\sigma}} / 3 $ taken to be positive in compression
\begin{equation}
\tensor{\sigma}= - p \, \id + \deviatoric{\tensor{\sigma}}
\; .
\end{equation}
%
For incompressible Newtonian fluids, the deviatoric part is customary expressed as
\begin{equation}
\deviatoric{\tensor{\sigma}}=2\mu\, \sym{\gradient{\vect{v}}}
\end{equation} 
with $\mu$ representing the fluid viscosity, here assumed to be isotropic and constant.
Momentum balance (\ref{eq:MomentumBalance}) thus becomes
\begin{equation}
-\gradient{ p }
+
\mu
\,
\divergence{ \gradient{ \vect{v} } }
+
\frac{\mu}{\rho}
\,
\gradient{s}
=
\vect{0}
\; .
\label{eq:NavierStokes}
\end{equation}
%
Eq. \eqref{eq:NavierStokes} is the usual Navier-Stokes equation neglecting inertia and body forces and accounting for the source mass contribution $s( {\vect{x}},t )$ in the current configuration, i.e. $ ({\vect{x}},t) \in \mathcal{P}(t) \times [0,T] $.

\subsubsection{Lubrication equation}
\label{subsubsec:lubrEq}

In the realm of hydraulic fracture, it is well established to make use of the lubrication assumption \cite{Hamrock2004} to simplify Navier-Stokes equation \eqref{eq:NavierStokes}.
%
Since in normal conditions the characteristic dimension of the fluid film along the fracture surface is significantly greater than its thickness,
lubrication theory moves from the assumption that %
the injection of mass $s( {\vect{x}},t )$ is uniform across the thickness of the crack. 
One thus defines the mean pressure $p_w$, density $\rho_w$, viscosity $\mu_w$ across the thickness and the mass supply per unit area $s_w$ on the fracture surface at point ${{\vect{y}}} \in
\Gamma_w$ as the integrals
\begin{align}
\nonumber
&
p_w ( {\vect{y}},t ) = \frac{1}{w ( {\vect{y}},t )} \int_{-w/2}^{w/2} \,  p( {\vect{x}},t ) \, {\rm d} x_3
\; ,
\qquad
\rho_w ( {\vect{y}},t ) = \frac{1}{w ( {\vect{y}},t )} \int_{-w/2}^{w/2} \,  \rho( {\vect{x}},t ) \, {\rm d} x_3
\; ,
\\ &
\mu_w ( {\vect{y}},t ) = \frac{1}{w ( {\vect{y}},t )} \int_{-w/2}^{w/2} \,  \mu( {\vect{x}},t ) \, {\rm d} x_3
\; ,
\qquad
s_w ( {\vect{y}},t ) =  \int_{-w/2}^{w/2} \,  s( {\vect{x}},t ) \, {\rm d} x_3
\; .
\label{eq:lubr_p_and_s}
\end{align}
The fluid flow is parallel to the crack plane and limiting the focus to fractures filled by Newtonian viscous fluids, 
the so called {\em{lubrication equation}} reads:
\begin{equation}
\frac{\rm d}{{\rm d}t}  w ( {\vect{y}},t ) 
=  
\frac{s_w ( {\vect{y}},t ) }{\rho_w ( {\vect{y}},t ) }
-
\surfacedivergence{
 \frac{w^3( {\vect{y}},t )}{12 \, \mu_w( {\vect{y}},t )}
  \left(
  -\surfacegradient{ p_w ( {\vect{y}},t )}{\Gamma_w}
  + \frac{ \mu_w }{ \rho_w } 
  \surfacegradient{  \frac{ s_w ( {\vect{y}},t )  }{ w ( {\vect{y}},t ) }  }{\Gamma_w}
  \right) 
 }{\Gamma_w}
 \; .
\label{eq:Lubrication1}
\end{equation}
%
In equation \eqref{eq:Lubrication1} we explicitly used the notion of surface divergence $\surfacedivergence{}{\Gamma_w}$ and gradient $\surfacegradient{}{\Gamma_w}$ to point out that the Navier-Stokes equation \eqref{eq:NavierStokes}, defined over a three-dimensional volume in the current configuration ${\vect{x}} \in \mathcal{P}(t)$,  
is now restricted to the fracture surface. We will not proceed with this notation henceforth, for the sake of readability we will make no distinction between surface and volumetric operators,
warning the reader at convenience.
%
The mathematics that lead to \eqref{eq:Lubrication1} have been summarized in appendix \ref{app:DimLubdEq}.

Focusing to the penny shaped fracture of radius $a(t)$ filled by a Newtonian viscous fluid depicted in Fig. \ref{Fig::fracture_draw},
we may switch to polar coordinates $0 \le \vartheta \le 2\pi$ , $0 \le r \le \ell(t)$ assuming the height of the crack to be small compared to its radius $w(r,t)\ll a(t)$.
Equation (\ref{eq:Lubrication1}) can be conveniently reshaped exploiting axis-symmetric conditions
\begin{eqnarray}
12 \mu_w
\; 
\frac{\rm d}{{\rm d}t}  w
=
12 \mu_w
\;
\frac{s_w}{\rho_w}
-
\frac{1}{r} \, 
\frac{\partial}{\partial r}
\left[ 
r w^3
\left(
-  \frac{\partial p_w}{\partial r }
+ 
\frac{\mu_w}{\rho_w}
\,
\frac{\partial }{\partial r} \, \frac{s_w}{w}
\right)
\right]
\; ,
\label{eq:Lubrication2}
\end{eqnarray}
%
with all fields estimated at $0 \le r \le \ell(t)$.

\subsection{Fronts tracking}
\label{subsec::FrontsTracking}

Inertia forces have not been considered in the analysis of the elastic body that encloses the fracture as well as of the fluid that flows in the fracture itself. Nonetheless, the hydro-fracture problem is time-dependent, because the fluid and fracture domains change with time. Accordingly, the problem formulation must include tracking equations for the fluid and for the crack fronts. They will be handled in the next two sections, making explicit reference to the penny shaped fracture of radius $a(t)$ filled by a Newtonian viscous fluid depicted in Fig. \ref{Fig::fracture_draw}.

\subsubsection{Fluid front tracking}
\label{subsec::FluidFrontTracking}

The driving force for the propagation of the penny shape fracture is an injection of fluid, which is modeled in this work as a mass smoothly supplied, in space and time, in a sub-region of the fracture. For the sake of symmetry, the fluid flux at the center of the circular crack must be vanishing. Impermeability condition implies that no fluid leak occurs across the lips of the fracture. Accordingly, the contribution along the boundary ${\partial\mathcal{P}(t)}$ in the mass balance (\ref{eq:MassBalance}) vanishes. 


Under the hypothesis of fluid incompressibility ($\partial \rho/\partial t=0$), Reynold's contribution of equation (\ref{eq:MassBalance}) reads
\begin{eqnarray}
\frac{ {\rm d} }{ {\rm d} t}
\int_{\mathcal{P}(t)} \rho \, {\rm d} v
&=&
\int_{\partial\mathcal{P}(t)}\rho \vect{v} \cdot\vect{n} \, {\rm d} a
\\
&=&
\nonumber
2\pi 
\,
\int_0^{\ell(t)}
\rho_w \, r \, \frac{ {\rm d} w}{ {\rm d} t} \, {\rm d} r
+
2\pi
\,
\int_0^{w(\ell(t),t)}
\rho_w \, \ell(t) \, \frac{\partial\ell(t)}{\partial t} \, {\rm d} r
\; .
\end{eqnarray}
%
%
The mass balance \eqref{eq:MassBalance}, therefore, yields the following estimation of the fluid front velocity
\begin{equation}
\label{eq:FluidFrontVelocity}
\frac{\partial\ell(t)}{\partial t} 
=
\frac{1}{\rho_w \,  \, \ell(t) \; w( \ell(t),t ) }
\int_0^{\ell(t)}
\left( 
  s_w - \rho_w \,  \frac{ {\rm d} w}{ {\rm d} t}
\right)
\; r 
\, {\rm d}r
\end{equation}
Integration of lubrication equation (\ref{eq:Lubrication2}) along the fluid film length after multiplication by $r$, leads to
\begin{equation}
\label{eq:IntLubrEq}
\int_0^{\ell(t)}
\rho_w \,  \frac{ {\rm d} w}{ {\rm d} t}
\; r
\, {\rm d} r
=
\int_0^{\ell(t)} 
\,
r \,
s_w
\,
{\rm d}r
-
\frac{1}{12\mu_w}
\left.
\left[
r  
\,
w^3
\left(
-  \frac{\partial p_w}{\partial r }
+ 
\frac{\mu_w}{\rho_w}
\,
\frac{\partial }{\partial r} \, \frac{s_w}{w}
\right)
\right]
\right|_0^{\ell(t)} 
\end{equation}
By plugging expression (\ref{eq:IntLubrEq}) into (\ref{eq:FluidFrontVelocity}), one finally has
\begin{equation}
\label{eq:FluidFrontVelocity2}
\frac{\partial\ell(t)}{\partial t}
=
\frac{ w^2( \ell(t), t )  }{12\mu_w}
\left(
\left.
 - 
\frac{\partial p_w }{\partial r}
\right|_{\ell(t)}
+
\left.
\frac{\mu_w}{\rho_w}
\,
\frac{\partial }{\partial r} \, \frac{s_w}{w}
\right|_{\ell(t)}
\right)
\; .
\end{equation}
%
%
%
\subsubsection{Crack front tracking}
\label{subsec::CrackFrontTracking}

The evolution of fractures in brittle or embrittled materials has been recently investigated in the theory of standard dissipative processes \cite{SalvadoriJMPS08,SalvadoriJMPS10}.
Variational formulations were stated \cite{ SalvadoriCariniIJSS2011,SalvadoriFantoni2013,
SalvadoriFantoniIJSS2014}, characterizing the crack front quasi-static velocity as the minimizer of constrained quadratic functionals.
The resulting implicit in time crack tracking algorithms \cite{SalvadoriFantoniJECS2014, SalvadoriFantoniJMPS2016} suffer from a major drawback, which limited the interest to their theoretical content and confined their computational relevance.
Crack tracking algorithms are formulated in terms of so-called {\itshape fundamental kernels} \cite{Lazarus2011}, whose values depend exclusively upon the shape of the crack front and whose analytical expression or even an accurate approximation is still missing for generic crack front shapes.
To overcome the issue, novel theoretical studies and resulting explicit in time crack tracking algorithms have been recently proposed in \cite{SalvadoriEtAlJMPS2019} for LEFM: 
%
they are grounded on a viscous regularization of the quasi-static fracture propagation problem and allow computing finite elongations of the crack front without fundamental kernels.
\\
\\
\textbf{Thermodynamic framework for fracture propagation}
\\
According to the Maximum Energy Release Rate criterion, the onset of fracture propagation at a point $s$ along the crack front at time $t$ can be mathematically expressed through the identity
\begin{equation}
\varphi(s,t)=G(s,t)-G_C = 0.
\label{eq:OnsetCrackProp}
\end{equation}  
Since energy dissipation during crack propagation is confined along the crack front, $G(s,t)$ in eq. (\ref{eq:OnsetCrackProp}) represents the energy release rate. For fractures that propagate in their own plane in pure mode I conditions, the energy release rate is related to the mode I SIF $K_I$ through Irwin's formula \cite{Irwin1958}, that for isotropic material reads
\begin{equation}
G(s,t)=\frac{ K_I^2(s,t)}{E'}.
\label{eq:irwin}
\end{equation}
$G_C$ in eq. (\ref{eq:OnsetCrackProp}) is the fracture energy, i.e. the dissipation per unit created crack surface here assumed as constant in time and space and related to the material fracture toughness $K_I^C$ through 
\begin{equation}
G_C=\frac{{K_I^C}^2}{E'}.
\label{eq:GC}
\end{equation} 
A safe equilibrium domain  can thus be defined, reminiscence of the elastic domain  in standard dissipative theory, expressing the condition for which cracks cannot advance. It is written as
\begin{equation}
\mathbb{E}=
\left\{
K_I(s,t)\in\mathbb{R}_0^+\left|\varphi(s,t)<0\right.\right\}
\label{eq:safeEquilibriumDomain}
\end{equation}
The boundary $\partial\mathbb{E}$ of $\mathbb{E}$ identifies the onset of crack propagation and it is a reminiscence of the yield surface in plasticity.
It is defined as
\begin{equation}
\partial\mathbb{E}=\left\{K_I(s,t)\in\mathbb{R}_0^+\left|\varphi(s,t)=0\right.\right\}
\label{eq:boundaryOfSafeEquilibriumDomain}
\end{equation}
SIFs shall belong to the closure of the safe equilibrium domain $\bar{\mathbb{E}}=\mathbb{E}\cup\partial\mathbb{E}$ and SIFs not satisfying this condition are ruled out.
The following inequalities,  dictated by the physics of the problem, govern fracture propagation under the assumption of irreversible crack growth
\begin{equation}
\varphi(s,t)\leq 0,
\hspace{0.5 cm}
\left.\frac{\partial a}{\partial t}\right|_{s,t}\geq 0,
\hspace{0.5 cm}
\varphi(s,t)\left.\frac{\partial a}{\partial t}\right|_{s,t}=0,
\label{eq:KarushKhunTucker_1}
\end{equation} 
where $\left.\frac{\partial a}{\partial t}\right|_{s,t}$ denotes the crack front quasi-static velocity at point $s$ and time $t$.
\\
\\
\textbf{ Fracture mechanics duality pairs}
\\
According to Griffith's theory, the power released at a point $s$ along the crack front as a consequence of fracture elongation is equal to the product between  the energy release rate $G(s,t)$ and the crack front  velocity at the same point $\left.\frac{\partial a}{\partial t}\right|_{s,t}$.
Such a power equals the product between  the SIF $K_I(s,t)$, playing the role of a thermodynamic force in the rigid plastic analogy, and the rate of its thermodynamically conjugated variable $v(s,t)$, namely
\begin{equation}
G(s,t)\left.\frac{\partial a}{\partial t}\right|_{s,t}=K_I(s,t)\,v(s,t).
\label{eq:DissipationEquality}
\end{equation}
A maximum dissipation principle can therefore be stated, as usual in standard dissipative systems theory. It holds
\\
{\itshape At any point $s$ along the crack front, for given $v(s,t)$, among all possible SIFs $k_I\in\bar{\mathbb{E}}$, the function}
\begin{equation}
\mathcal{D}\left(k_I;v\right)=k_I\,v
\label{eq:Constrained_D_a}
\end{equation}
{\itshape attains its maximum  for the actual SIFs $K_I$:}
\begin{equation}
\mathcal{D}\left(K_I;v\right)=max_{k_I\in\bar{\mathbb{E}}}\mathcal{D}\left(k_I;v\right).
\label{eq:Constrained_D}
\end{equation}
Maximality condition (\ref{eq:Constrained_D}) implies associativity
\begin{equation}
v(s,t)=\left.\frac{\partial\varphi}{\partial K_I}\right|_{s,t}\dot{\lambda}(s,t),
\label{eq:ConstrainedFlowRule}
\end{equation}
and Karush-Khun-Tucker  loading/unloading conditions
\begin{equation}
\varphi(s,t)\leq 0,
\hspace{0.5 cm}
\dot{\lambda}(s,t)\geq 0,
\hspace{0.5 cm}
\varphi(s,t)\,\dot{\lambda}(s,t)=0,
\label{eq:KarushKhunTucker_2}
\end{equation}
which are equivalent to inequalities (\ref{eq:KarushKhunTucker_1}).
If the onset of crack propagation is equivalently written in the form $\varphi(s,t)=K_I(s,t)-K_I^C$, from normality condition (\ref{eq:ConstrainedFlowRule}) one obtains that $v(s,t)$ equals the plastic multiplier $\dot{\lambda}(s,t)$.
Thermodynamic consistency
\begin{equation}
K_I(s,t)v(s,t)\geq 0
\label{eq:ClDuInequality}
\end{equation} 
in therefore satisfied in view of positive definiteness of $K_I$ and complementarity laws (\ref{eq:KarushKhunTucker_2}).
Plugging eq. (\ref{eq:irwin}) into (\ref{eq:DissipationEquality}), one has
\begin{equation}
v(s,t)=\dot{\lambda}(s,t)=\frac{K_I}{E'}\left.\frac{\partial a}{\partial t}\right|_{s,t},
\label{eq:crackFrontVelocity}
\end{equation}
thus relating $v(s,t)$ to the crack front quasi-static velocity $\left.\frac{\partial a}{\partial t}\right|_{s,t}$ at point $s$ of the crack front.
\\
\\
\textbf{Viscous regularization of the quasi-static fracture propagation problem}
\\
As usual in standard dissipative systems, viscosity can be interpreted as a regularization of rate-independent formulations \cite{SIMOHUGHES}.  
Classical viscous constitutive models, in fact, can emanate from optimality conditions imposed on a regularized penalty form of the maximum dissipation function.
The viscosity we will account for, hereafter denoted with $\eta$, is not a {\em{material viscosity}}, rather it is a numerical regularizing parameter that solely enters the algorithm in predicting the geometrical evolution of the fracture front.
In other words, $\eta$ does neither affect the overall elastic behavior of the host material nor the order of singularity of the stress field ahead of the crack.
In view of the viscous regularization, Karush-Khun-Tucker conditions (\ref{eq:KarushKhunTucker_1})  do not hold anymore and $K_I \not\in \bar{\mathbb{E}}$ are admissible.
Denoting viscous SIFs with $K_I^\eta$ and with $v^\eta(s,t)$ the rate of their conjugate internal variable, dissipation identity (\ref{eq:DissipationEquality}) transforms into
\begin{equation}
G(s,t)\left.\frac{\partial a }{\partial t}\right|_{s,t}=K_I^\eta(s,t)\,v^\eta(s,t).
\label{eq:DissipationEquality_2}
\end{equation}
An unconstrained formulation of the maximum dissipation principle thus applies, reminiscent, for example, of Perzina-type visco-plastic model.
It reads
\\

{\itshape At any point $s$ along the crack front, for given viscosity $\eta$ and $v^\eta(s,t)$, among all possible SIFs $k_I\in\mathbb{R}_0^+$, the function} 
\begin{equation}
\mathcal{D}\left(k_I;v^\eta\right)=k_I\,v^\eta(s,t)-\frac{\varphi^2(s,t)}{2\eta}\mathcal{H}[\varphi],
\label{eq:UNConstrained_D}
\end{equation}
{\itshape attains its maximum for the actual SIF $K_I^\eta$:}
\begin{equation}
\mathcal{D}\left(K_I^\eta;v^\eta\right)=max_{k_I\in\mathbb{R}_0^+}\mathcal{D}\left(k_I;v^\eta\right),
\label{eq:UNConstrainedDissipationEquality}
\end{equation}
where $\mathcal{H}[\varphi]$ is the Heaviside function that holds 1 if $\varphi>0$ and vanishes otherwise, and  the numerical viscosity $\eta\in(0,+\infty)$ represents the penalty parameter of the unconstrained formulation.
At this point, regularized form of normality rule (\ref{eq:ConstrainedFlowRule}) is obtained enforcing the vanishing of the Gateaux derivative of functional (\ref{eq:UNConstrained_D}):
\begin{equation}
v^\eta(s,t)=\left.\frac{\partial\varphi}{\partial K_I^\eta}\right|_{s,t}\frac{\varphi(s,t)}{\eta}\mathcal{H}[\varphi].
\label{eq:UNConstrained_Flow_rule}
\end{equation}
Crack front quasi-static velocity amounts at
\begin{equation}
\left.\frac{\partial a}{\partial t}\right|_{s,t}=\frac{E'}{\eta}
\left(
1-\frac{K_I^C}{K_I^\eta}
\right)
\mathcal{H}
\left[
K_I^\eta-K_I^C
\right]
\, ,
\label{eq:CrackFrontVelocity}
\end{equation}
in view of identity (\ref{eq:DissipationEquality_2}).
%
%
It is a classical argument in optimization theory that, provided sufficient smoothness, $K_I^\eta\rightarrow K_I$ as $\eta\rightarrow 0$, thus providing the convergence of unconstrained formulation (\ref{eq:UNConstrained_D}) to the solution of the constrained counterpart (\ref{eq:Constrained_D_a}). For this sake, the apex $^\eta$ will be omitted henceforth.
Performed regularization of the quasi-static fracture propagation problem allows obtaining an effective expression  for the quasi-static crack front velocity materialized in eq. (\ref{eq:CrackFrontVelocity}), which can be implemented in crack tracking algorithms by imposing Griffith's criterion at every propagation step.
The approximation \eqref{eq:CrackFrontVelocity} of the crack front velocity relies on the ``distance'' between the stress intensity factor $K_I(a(t))$ and the critical one, $K_I^C$. Since the latter cannot be overcome, crack tracking algorithms aim at finding $a(t)$ such that $K_I(a(t)) =K_I^C$, which is an equilibrium condition for $a(t)$ in view of eq. \eqref{eq:CrackFrontVelocity}. Accordingly, the knowledge of $K_I(a(t))$ is required at any length $a$. Numerical approximation of the SIF can be achieved either by finite elements \cite{ZammarchiEtAlCMAME2017} moving from the weak form \eqref{WeakNavier}, or via boundary elements \cite{SalvadoriGray06} exploiting asymptotic analysis for \eqref{Overview:IntEq}. In simple cases, as for pressurized penny shaped cracks under axis-symmetric conditions, analytical solutions are available \cite{Sih1986}, too:
\begin{equation}
\label{eq:SIFKI}
K_I(a(t))=\frac{2}{\pi\sqrt{a(t)}}
\int_0^{a(t)}
\frac{r \, p_w (r,t)}{\sqrt{a^2(t)-r^2}}
\,
{\rm d} r
\; .
\end{equation}
%

%
%
\section{Space and time discretization of governing equations}
\label{Sec::Discretization}

%

The elastic response \eqref{eq:wNostra} and the lubrication equation \eqref{eq:Lubrication2} provide the crack opening $w$ and the pressure distribution $p_w$ at equilibrium, complemented by equations \eqref{eq:FluidFrontVelocity} or \eqref{eq:FluidFrontVelocity2} and by eq. \eqref{eq:CrackFrontVelocity}.
They establish  a complete system of equations to determine the problem unknown fields $w(r,t)$, $p_w(r,t)$, $\ell(t)$, and $a(t)$ along the penny shaped crack surface.

The axis-symmetric weak form of the lubrication equation can be obtained by multiplication of (\ref{eq:Lubrication2}) for a time-independent test function $\hat{p}(r)$ 
and by integration over the fracturing fluid plane, thus obtaining
\begin{align}
\label{eq:weakform2}
12\mu_w 
\;
\int_0^{\ell(t)} 
\, \hat{p}(r) \, 
\frac{\rm d}{{\rm d}t}  w(r,t) \,
r \,
{\rm d} r
=
\;
&
12 \,
\frac{\mu_w}{\rho_w}
\, 
\int_0^{\ell(t)} 
\, \hat{p}(r) \, s_w(r,t) \, r \,
{\rm d} r
+
\nonumber\\
&
-\int_0^{\ell(t)}
\frac{\partial}{\partial r}
\left\{
 r w^3(r,t)
\left(
-  \frac{\partial p_w(r,t)}{\partial r }
+ 
\frac{\mu_w}{\rho_w}
\,
\frac{\partial }{\partial r} \, \frac{s_w(r,t)}{w}
\right)
\hat{p}(r)
\right\}
\,
{\rm d} r
+
\nonumber\\
&
+
\int_0^{\ell(t)} r w^3(r,t)
\left(
-  \frac{\partial p_w(r,t)}{\partial r }
+ 
\frac{\mu_w}{\rho_w}
\,
\frac{\partial }{\partial r} \, \frac{s_w(r,t)}{w}
\right)
\frac{\partial\hat{p}(r)}{\partial r} \,
{\rm d} r
\end{align}
Since $\hat{p}(\ell)=0$, it holds
\begin{equation}
\left.
r w^3(r,t)
\left(
-  \frac{\partial p_w(r,t)}{\partial r }
+ 
\frac{\mu_w}{\rho_w}
\,
\frac{\partial }{\partial r} \, \frac{s_w(r,t)}{w}
\right)
\hat{p}(r)
\right|_0^{\ell(t)}
=
0
\; .
\end{equation}
The weak form of lubrication equation eventually writes as:
\begin{align}
\label{eq:weakform3}
12\mu_w 
\;
\int_0^{\ell(t)} 
\, \hat{p}(r) \, 
\frac{\rm d}{{\rm d}t}  w(r,t) \,
r \,
{\rm d} r
=
&
\;
12 \,
\frac{\mu_w}{\rho_w}
\, 
\int_0^{\ell(t)} 
\, \hat{p}(r) \, s_w(r,t) 
\, r \,
{\rm d} r
+
\nonumber\\
&
-
\int_0^{\ell(t)}
\frac{\partial \hat{p}(r)}{\partial r}
\;  w^3(r,t) \;
\frac{\partial p_w(r,t)}{\partial r}
\,
r \,
{\rm d} r
+
\nonumber\\
&
+
\frac{\mu_w}{\rho_w}
\,
\int_0^{\ell(t)}  
\frac{\partial \hat{p}(r)}{\partial r}
\; w^2(r,t) \;
\frac{\partial s_w(r,t)}{\partial r}
\, r \,
{\rm d} r
\nonumber\\
&
-
\int_0^{\ell(t)}  
\frac{\partial \hat{p}(r)}{\partial r}
\; w(r,t) \, s_w(r,t) \;
\frac{\partial w(r,t)}{\partial r}
\, r \,
{\rm d} r
\; .
\end{align}
The weak form  \eqref{eq:weakform3} naturally leads to a semi-discrete problem, by approximating the unknown fields as a product of separated variables, by means of a basis $ \{ \psi_k (r) \}$  of spatial shape functions and nodal unknowns that depend solely on time. For the sake of readability, the apex $_w$ will not be written any further.
One has
\begin{equation}
\label{eq:pressurefieldDiscretization}
p(r,t)=\sum_{k=1}^{N_{nod}}
\psi_k(r) \, p_k(t),
\end{equation}
with $N_{nod}$ the number of nodal unknowns $p_k(t)$.
%
%
The crack opening $w(r,t)$ can be discretized by means of integral (\ref{eq:wNostra}), i.e.
\begin{equation}
\label{eq:wDisctretization}
w(r, t)=\frac{8}{\pi E'} \, \sum_{k=1}^{N_{nod}} \, A_k(a(t),\ell(t);r,t) \; p_k(t),
\end{equation}
where discrete integral influence functions 
\begin{equation}
\label{eq:InfluenceFunctionAk}
 A_k( a(t),\ell(t);r,t )
 =
 \int_r^{a(t)} 
 \frac{1}{\sqrt{z^2-r^2}}
 \int_0^z
 \frac{y \, \psi_k(y)}{\sqrt{z^2-y^2}}
 \,
 {\rm d}y
 \,
 {\rm d}z,
\end{equation}
relate the crack opening to the nodal values of pressure and depend upon the crack and the fluid lengths at time $t$. Their computation is detailed in Appendix \ref{app:Ak} for linear shape functions.

\bigskip
In view of (\ref{eq:wDisctretization}), the discrete pressure field \eqref{eq:pressurefieldDiscretization} governs the whole problem, through time-dependent nodal unknowns $p_k(t)$.
The total derivative of crack opening ${{\rm d} w} / { {\rm d}t} $ that appears in the lubrication equation (\ref{eq:Lubrication2})  can be discretized as
\begin{equation}
\frac{ {\rm d} w(r,t)}{  {\rm d} t}
=
\frac{8}{\pi E'} \sum_{k=1}^{N_{nod}}
A_k(r,t) \, \frac{\partial p_k(t)}{\partial t} +
\frac{\partial A_k(r,t)}{\partial a} \; \frac{\partial a(t)}{\partial t} \, p_k +
\frac{\partial A_k(r,t)}{\partial \ell} \; \frac{\partial \ell(t)}{\partial t} \, p_k
\end{equation}
where dependence of $A_k$ upon $a(t)$ and $ \ell(t)$ will be omitted henceforth for the sake of brevity.
Discrete weak form (\ref{eq:weakform3}) eventually reads
\begin{eqnarray}
\label{eq:weakform4}
&&
\frac{96 \mu }{\pi E'}\sum_{k=1}^{N_{nod}}
\left[
 \int_0^{\ell(t)}
 \, \psi_i(r) \, 
 A_k(r,t) 
  \, r \, {\rm d} r
  \; 
\frac{\partial p_k(t)}{\partial t}
+
\int_0^{\ell(t)}
 \, \psi_i(r) \, 
\frac{\partial A_k(r,t)}{\partial a}
\frac{\partial a(t)}{\partial t}
  \, r \, {\rm d} r
  \;
 p_k(t)
 +
 \right.
 \nonumber\\
 &&
 \hspace{1.5 cm}
 \left.
 \int_0^{\ell(t)}
 \, \psi_i(r) \, 
\frac{\partial A_k(r,t)}{\partial \ell}
\frac{\partial \ell(t)}{\partial t}
  \, r \, {\rm d} r
  \;
 p_k(t)
 \right]
 +
 \nonumber\\
 &&
  +
 \left(
 \frac{8}{\pi E'} 
 \right)^3 
 \,
 \sum_{k=1}^{N_{nod}}
 \left[
 \int_0^{\ell(t)}
 \,
 \frac{\partial \psi_i(r)}{\partial r} 
 \;
 \left(
 \sum_{h=1}^{N_{nod}}
 A_h(r,t)p_h(t)
 \right)^3 
 \;
 \frac{\partial \psi_k(r)}{\partial r}
  \, r \, {\rm d} r
  \;
 p_k(t)
 \right]
 +
 \nonumber\\
 &&
 -
  \left(
 \frac{8}{\pi E'} 
 \right)^2
 \,
 \frac{\mu}{\rho}
 \int_0^{\ell(t)}
 \frac{\partial \psi_i(r)}{\partial r} 
 \,
 \left(
 \sum_{h=1}^{N_{nod}}
  A_h(r,t)p_h(t)
 \right)^2
 \,
 \frac{\partial s_w(r,t)}{\partial r}
  \, r \, {\rm d} r
  \;
 +
  \nonumber\\
 &&
 +
   \left(
 \frac{8}{\pi E'} 
 \right)^2
 \,
 \int_0^{\ell(t)}
 \frac{\partial \psi_i(r)}{\partial r} 
 \,
 \left(
 \sum_{h=1}^{N_{nod}}
 A_h(r,t)p_h(t)
 \right)
 s_w(r,t)
 \left(
 \sum_{j=1}^{N_{nod}}
\frac{\partial A_j(r,t) }{\partial r} 
  p_j(t)
 \right)
  \, r \, {\rm d} r
  \;
 \nonumber\\
 && 
 =
 \frac{12\mu}{\rho} 
 \int_0^{\ell(t)}
 \, \psi_i(r) \, 
 s_w(r,t)
  \, r \, {\rm d} r
  \;
  \qquad
 \hspace{0.5 cm}
 i=1,...,N_{nod}
 \nonumber\\
\end{eqnarray}
Computation of  $\partial A_k /\partial a $, $\partial A_k /\partial \ell  $ and $\partial A_h /\partial r $ in equation (\ref{eq:weakform4}) is detailed in Appendix \ref{app:Ak} for linear shape functions.

\bigskip
In order to obtain a full discretization of the weak form \eqref{eq:weakform4}, define $t_n = n \, \Delta t$ with $n=0,1,...$, with $\Delta t > 0$ being the time step. 
The differential equation (\ref{eq:weakform4}) in time can be solved via Backward Euler Method, leading to solution at time $t_{n+1}$ computed through the following equation

\begin{eqnarray}
\label{eq:BackEul}
&&
\nonumber
\frac{96 \mu }{\pi E'}\sum_{k=1}^{N_{nod}}
\left[
 \int_0^{\ell(t_{n+1})}
 \, \psi_i(r) \, 
 A_k(r,t_{n+1}) 
  \, r \, {\rm d} r
  \; 
\frac{p_k(t_{n+1})-p_k(t_n)}{\Delta t}
+
 \right.
\\
 &&
 \hspace{1.5 cm}
 \left.
\int_0^{\ell(t_{n+1})}
 \, \psi_i(r) \, 
\frac{\partial A_k(r,t_{n+1})}{\partial a}
\frac{\partial a(t_{n+1})}{\partial t}
  \, r \, {\rm d} r
  \;
 p_k(t_{n+1})
 +
 \right.
 \nonumber\\
 &&
 \hspace{1.5 cm}
 \left.
 \int_0^{\ell(t_{n+1})}
 \, \psi_i(r) \, 
\frac{\partial A_k(r,t_{n+1})}{\partial \ell}
\frac{\partial \ell(t_{n+1})}{\partial t}
  \, r \, {\rm d} r
  \;
 p_k(t_{n+1})
 \right]
 +
 \nonumber\\
 &&
  +
 \left(
 \frac{8}{\pi E'} 
 \right)^3 
 \,
 \sum_{k=1}^{N_{nod}}
 \left[
 \int_0^{\ell(t_{n+1})}
 \,
 \frac{\partial \psi_i(r)}{\partial r} 
 \;
 \left(
 \sum_{h=1}^{N_{nod}}
 A_h(r,t_{n+1})p_h(t_{n+1})
 \right)^3 
 \;
 \frac{\partial \psi_k(r)}{\partial r}
  \, r \, {\rm d} r
  \;
 p_k(t_{n+1})
 \right]
 +
 \nonumber\\
 &&
 -
  \left(
 \frac{8}{\pi E'} 
 \right)^2
 \,
 \frac{\mu}{\rho}
 \int_0^{\ell(t_{n+1})}
 \frac{\partial \psi_i(r)}{\partial r} 
 \,
 \left(
 \sum_{h=1}^{N_{nod}}
  A_h(r,t_{n+1})p_h(t_{n+1})
 \right)^2
 \,
 \frac{\partial s_w(r,t_{n+1})}{\partial r}
  \, r \, {\rm d} r
  \;
 +
  \nonumber\\
 &&
 +
   \left(
 \frac{8}{\pi E'} 
 \right)^2
 \,
 \int_0^{\ell(t_{n+1})}
 \frac{\partial \psi_i(r)}{\partial r} 
 \,
 \left(
 \sum_{h=1}^{N_{nod}}
 A_h(r,t_{n+1})p_h(t_{n+1})
 \right)
 s_w(r,t_{n+1})
 \left(
 \sum_{j=1}^{N_{nod}}
\frac{\partial A_j(r,t_{n+1}) }{\partial r} 
  p_j(t_{n+1})
 \right)
  \, r \, {\rm d} r
  \;
 \nonumber\\
 && 
 =
 \frac{12\mu}{\rho} 
 \int_0^{\ell(t_{n+1})}
 \, \psi_i(r) \, 
 s_w(r,t_{n+1})
  \, r \, {\rm d} r
  \;
  \qquad
 \hspace{0.5 cm}
 i=1,...,N_{nod}
\end{eqnarray}
%
%
%
with the crack front velocity expressed as
\begin{equation}
\frac{\partial a(t_{n+1})}{\partial t}=
\frac{E'}{\eta}
\left(
1-\frac{K_I^C}{K_I(a(t_{n+1}))}
\right)
\mathcal{H}
\left[
K_I(a(t_{n+1}))-K_I^C
\right]
%
\end{equation}
and the fluid front velocity in the form \eqref{eq:FluidFrontVelocity2}.
%
The four unknowns ${p}(r,t_{n+1}),{w}(r,t_{n+1}),a(t_{n+1})$, and $\ell(t_{n+1})$ of non the linear problem (\ref{eq:BackEul}) have been determined by means of a staggered approach, as detailed in box \ref{box:algobenchmark}. 
Therein, $\mu$ and $j$ stand for the counter of the outer loop devoted to compute the crack radius and the counter of the inner loop devoted to compute the fluid front position, respectively.

\begin{tcolorbox}[adjusted title= Box1: Meta-code for the  propagation of a penny shaped hydraulic fracture with fluid lag]
For $n=1:Ntimesteps$
\\
\\
$ \hspace{0.5cm}\mbox{}$
\,\,$\bullet\,$trigger with $\mu=0: a_{[0]}^{(n+1)}=a^{(n)}$
\\
\\
$ \hspace{0.5cm}\mbox{}$
\,\,\,\,$\bullet\,$trigger with $j=0:\ell^{(n+1)}_{[0]}=\ell^{(n)},\vect{p}_{[0]}^{(n+1)}=\vect{p}^{(n)}$
\\
\\
$ \hspace{0.5cm}\mbox{}$
\,\,\,\,\,\,\,\,\,\,
$\bullet\,$solve for $\ell_{[j]}^{(n+1)}$ at $a_{[\mu-1]}^{(n+1)}$:
\\
\\
$ \hspace{0.5cm}\mbox{}$
\,\,\,\,\,\,\,\,\,\,
$
\ell_{[j]}^{(n+1)}=\ell^{(n)}+\frac{\Delta t}{\ell_{[j-1]}^{(n+1)}\rho w^{(n+1)}(\ell_{[j-1]}^{(n+1)})}\int_0^{\ell_{[j-1]}^{(n+1)}}\left(s_w^{(n+1)}(r)-\rho\frac{w_{[j-1]}^{(n+1)}-w^{(n)}}{\Delta t}\right)r\,dr
$
\\
\\
$ \hspace{0.5cm}\mbox{}$
\,\,\,\,\,\,\,\,\,\,
$\bullet\,$solve for $\vect{p}_{[j]}^{(n+1)}$:
\\
\\
$ \hspace{0.5cm}\mbox{}$
\,\,\,\,\,\,\,\,\,\,
$A_{k[\mu-1]}^{(n+1)}(r)=\int_r^{a_{[\mu-1]}^{(n+1)}}\frac{1}{\sqrt{z^2-r^2}}\int_0^z\frac{y\psi_k(y)}{\sqrt{z^2-y^2}}\,dy\,dz,$
\\
\\
$ \hspace{0.5cm}\mbox{}$
\,\,\,\,\,\,\,\,\,\,
$\frac{96\mu}{\pi E'}\sum_{k=1}^{N_{nod}}\left[\int_0^{\ell_{[j]}^{(n+1)}}\psi_i(r)A_{k[\mu-1]}^{(n+1)}(r)r\,dr\frac{p_{k[j]}^{(n+1)}-p_k^{(n)}}{\Delta t}\right.+$
\\
\\
$ \hspace{0.5cm}\mbox{}$
\,\,\,\,\,\,\,\,\,\,
$\int_0^{\ell_{[j]}^{(n+1)}}\psi_i(r)\frac{ \partial  A_{k[\mu-1]}^{(n+1)}(r)}{\partial a}r\,dr\frac{a_{[\mu-1]}^{(n+1)}-a^{(n)}}{\Delta t}p_{k[j]}^{(n+1)}+$
\\
\\
$ \hspace{0.5cm}\mbox{}$
\,\,\,\,\,\,\,\,\,\,
$\left.\int_0^{\ell_{[j]}^{(n+1)}}\psi_i(r)\frac{ \partial  A_{k[\mu-1]}^{(n+1)}(r)}{\partial \ell}r\,dr\frac{\ell_{[j]}^{(n+1)}-\ell^{(n)}}{\Delta t}p_{k[j]}^{(n+1)}\right]=$
\\
\\
$ \hspace{1.5cm}\mbox{}$
\,\,\,\,\,\,\,\,\,\,
$=\frac{12\mu}{\rho}\int_0^{\ell_{[j]}^{(n+1)}}\psi_i(r) s_w^{(n+1)}(r)r\,dr+$
\\
\\
$ \hspace{1.5cm}\mbox{}$
\,\,\,\,\,\,\,\,\,\,
$
-\left(\frac{8}{\pi E'}\right)^3
\sum_{k=1}^{N_{nod}}\int_0^{\ell_{[j]}^{(n+1)}}
\frac{\partial \psi_i(r)}{\partial r}
\left(\sum_{h=1}^{N_{nod}} A_{h[\mu-1]}^{(n+1)}(r)p_{h[j-1]}^{(n+1)}\right)^3
\frac{\partial \psi_k(r)}{\partial r}
r\,dr p_{k[j]}^{(n+1)}$
\\
\\
$ \hspace{1.5cm}\mbox{}$
\,\,\,\,\,\,\,\,\,\,
$
+\left(\frac{8}{\pi E'}\right)^2
\frac{\mu}{\rho}
\int_0^{\ell_{[j]}^{(n+1)}}
\frac{\partial \psi_i(r)}{\partial r}
\left(\sum_{h=1}^{N_{nod}} A_{h[\mu-1]}^{(n+1)}(r)p_{h[j-1]}^{(n+1)}\right)^2
\frac{\partial s_w^{(n+1)}(r)}{\partial r}
r\,dr$
\\
\\
$ \hspace{0.5cm}\mbox{}$
\,\,\,\,\,\,\,\,\,\,
$
-\left(\frac{8}{\pi E'}\right)^2
\int_0^{\ell_{[j]}^{(n+1)}}
\frac{\partial \psi_i(r)}{\partial r}
\left(\sum_{h=1}^{N_{nod}}A_{h[\mu-1]}^{(n+1)}(r)p_{h[j-1]}^{(n+1)}\right)
 s_w^{(n+1)}(r)
 \left(\sum_{h=1}^{N_{nod}}\frac{\partial A_{h[\mu-1]}^{(n+1)}(r)}{\partial r}p_{h[j-1]}^{(n+1)}\right)
r\,dr$
\\
\\
$ \hspace{0.5cm}\mbox{}$
\,\,\,\,\,\,\,\,\,\,
$\forall\hspace{0.1cm}i=1,2,...,N_{nod} $
\\
\\
$ \hspace{1.5cm}\mbox{}$
\,\,\,\,\,\,\,\,\,\,
if $\frac{\left|\ell^{(n+1)}_{[j]}-\ell^{(n+1)}_{[j-1]}\right|}{\left|\ell^{(n+1)}_{[j]}\right|}<tol_{\ell}$
\\
\\
$ \hspace{1.5cm}\mbox{}$
\,\,\,\,\,\,\,\,\,\,
break
\\
\\
$ \hspace{1.5cm}\mbox{}$
\,\,\,\,\,\,\,\,\,\,
else $j++$
\\
\\
$ \hspace{1.5cm}\mbox{}$
\,\,
$\bullet\,$evaluate $K_{I[\mu]}^{(n+1)}=\frac{2}{\pi\sqrt{a^{(n+1)}_{[\mu-1]}}}\sum_{k=1}^{N_{nod}}\int_0^{a^{(n+1)}_{[\mu-1]}}
\frac{r\psi_k(r)}
{\sqrt{a^{(n+1)2}_{[\mu-1]}-r^2}}\,dr 
p_{k[j]}^{(n+1)}
$
\\
\\
$ \hspace{1.5cm}\mbox{}$
\,\,
$\bullet\,$update $a^{(n+1)}_{[\mu]}$:
\\
\\
$ \hspace{1.5cm}\mbox{}$
\,\,
$a^{(n+1)}_{[\mu]}=a^{(n+1)}_{[\mu-1]}+
\frac{\Delta t E'}{\eta}
\left(1-
\frac{K_I^C}{K_{I[\mu]}^{(n+1)}}
\right)
\mathcal{H}
\left[K_{I[\mu]}^{(n+1)}-K_I^C\right]$
\\
\\
$ \hspace{1.5cm}\mbox{}$
\,\,
$\bullet\,$if $\frac{\left|a^{(n+1)}_{[\mu]}-a^{(n+1)}_{[\mu-1]}\right|}{\left|a^{(n+1)}_{[\mu]}\right|}<tol_a$
\\
\\
$ \hspace{1.5cm}\mbox{}$
\,\,
break
\\
\\
$ \hspace{1.5cm}\mbox{}$
\,\,
else $\mu++$
\\
\\
end

\label{box:algobenchmark}
\end{tcolorbox}


\section{Model validation on a numerical benchmark.}
\label{Sec::Benchmark}

Consider a penny-shape crack with an initial radius $a^{(0)}=20\,mm$, immersed in a medium of density $\rho=10^3\,kg/m^3$, plane strain Young modulus $E'=4\cdot10^4 N/mm^2$ and fracture toughness $K_I^C= 0.5\,MPa\sqrt{m}$.
The crack hosts a viscous Newtonian fluid ($\mu=0.89\cdot 10^{-9}MPa\,s$), which initially fills a circular region of radius $\ell^{(0)}=10\,mm$ with an initial pressure $p_0=1\,MPa$. 
The volume of fluid trapped in the crack is enlarged by the following mass supply per unit area:
 \begin{equation*}
 s_w(r,t)=\frac{0.03}{2\pi}\textrm{exp}\left(-\frac{r^2}{32}\right)\left(1-\textrm{exp}\left(-2\,t\right)\right)
 \; .
 \end{equation*} 
The far field stress $\sigma_0$ has been considered as vanishing.
%
%
%
%
In this axis-symmetric case, mode I SIF has been computed integrating numerically eq. (\ref{eq:SIFKI}).

Simulations based upon the numerical scheme described in the former sections, with numerical viscosity $\eta$ set\footnote{The influence of such a parameter upon the performances of the crack tracking algorithm is discussed at large in \cite{SalvadoriEtAlJMPS2019}: 
the lower the value of $\eta$, the lower the number of iteration requested to satisfy a Griffith condition ($G=G_C$) at the crack front.
} to $0.1\,Nm^{-3}s$,  give rise to the pressure and opening fields depicted below. It is worth mentioning that terms containing $\frac{\partial A_k(r,t)}{\partial \ell}$ and $\frac{\partial A_k(r,t)}{\partial r}$  in equation (\ref{eq:BackEul}) result to be numerically negligible with respect to the other ones in the system resolution.

\begin{figure}[h]
\begin{subfigure}{0.33\textwidth}
  \includegraphics[width=.9\linewidth]{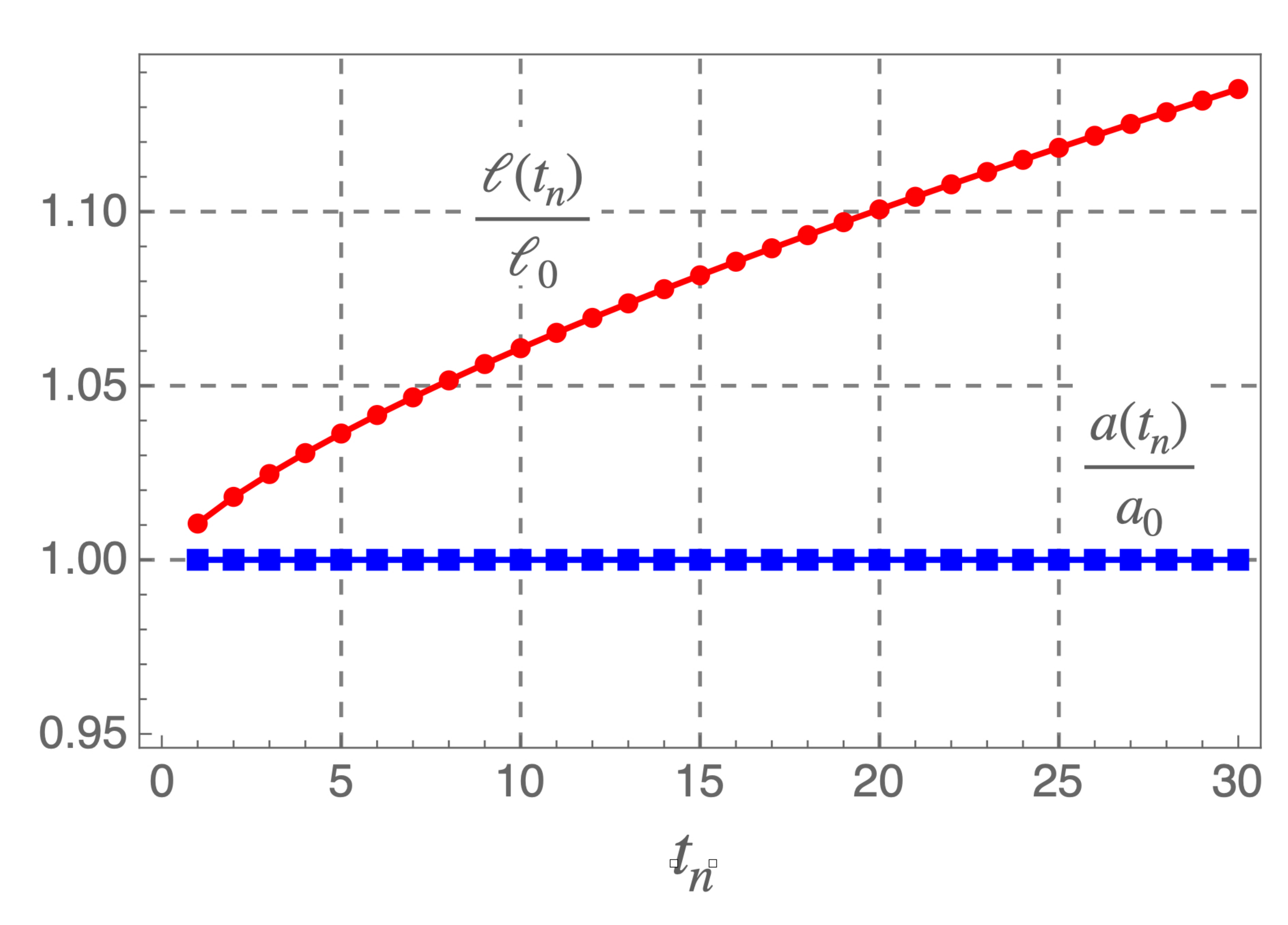}
\caption{}
\end{subfigure}
\begin{subfigure}{0.33\textwidth}
\centering
  \includegraphics[width=.9\linewidth]{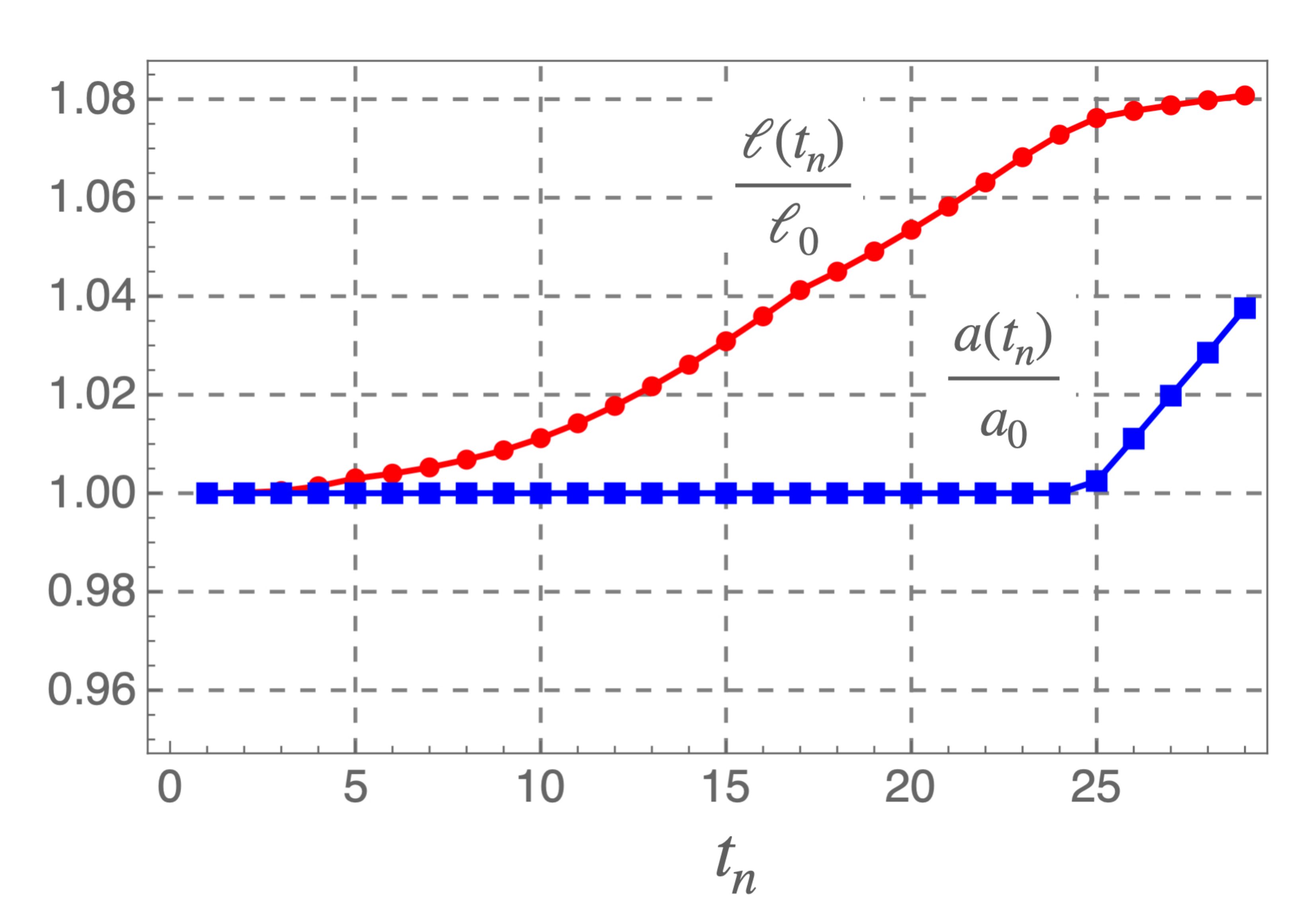}
\caption{}
\end{subfigure}
\begin{subfigure}{0.33\textwidth}
  \includegraphics[width=.9\linewidth]{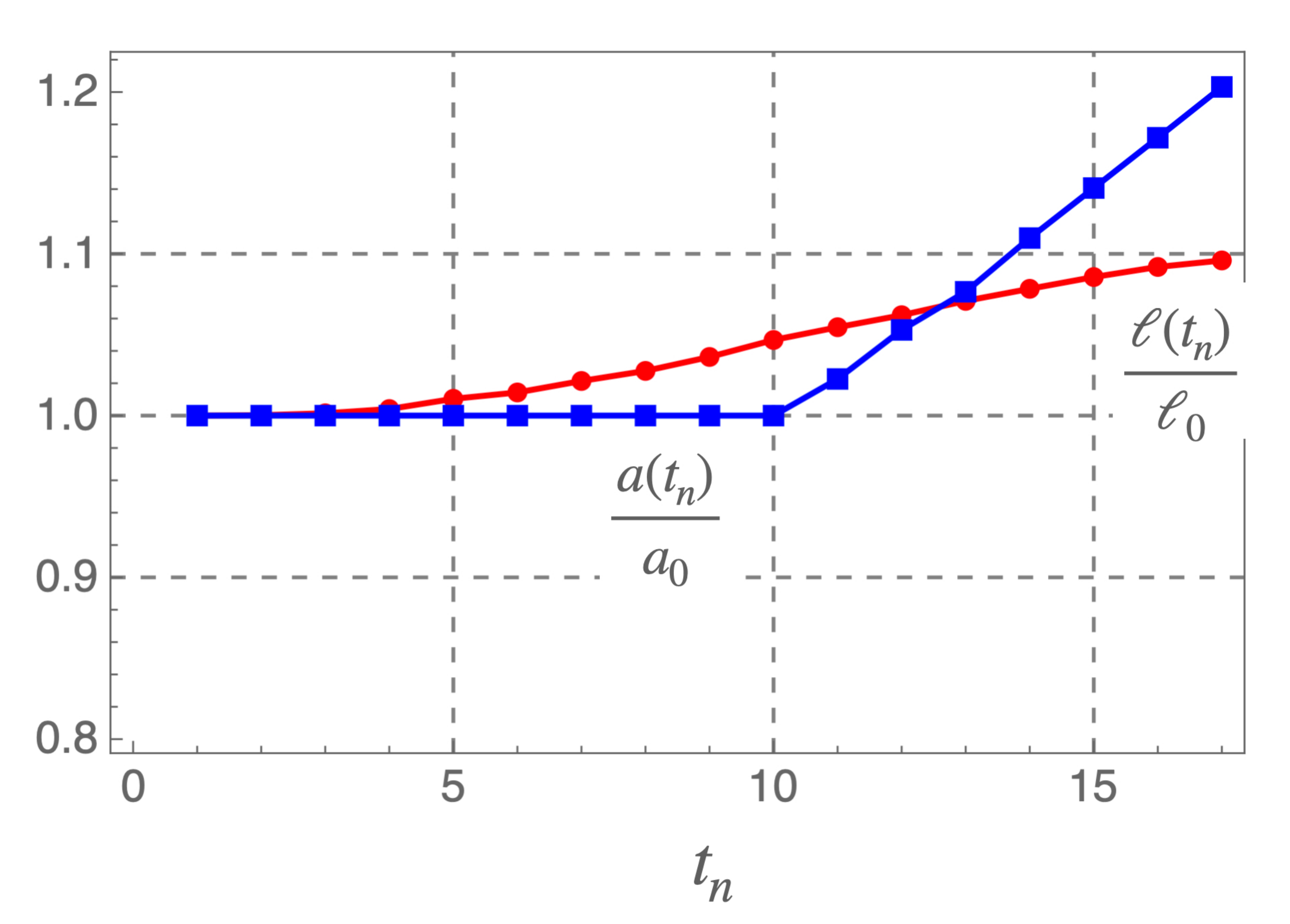}
\caption{}
\end{subfigure}
  \caption{\it Normalized fluid front $\ell^{(n)}$ (red curves) and crack  $a^{(n)}$ (blue curves) lengths vs time $t_n$. (a)$E'=4\cdot 10^2\, MPa$;  (b)  $E'=4\cdot 10^4\, MPa$; (c) $E'=4\cdot 10^6\, MPa$.}
  \label{Fig::ell-a}
\end{figure}

Figure \ref{Fig::ell-a}-(a) depicts the fluid ($\ell^{(n)}$) and the crack ($a^{(n)}$) lengths, normalized by the corresponding initial one ($n=0$), at several time steps for three surrounding media with increasing stiffness. In compliant media, the crack does not immediately advance whereas the fluid does so. Stiff surrounding media show the opposite behavior, with crack immediately growing and increasing fluid lag.

\begin{figure}[h]
\begin{subfigure}{0.5\textwidth}
\centering
  \includegraphics[width=.9\linewidth]{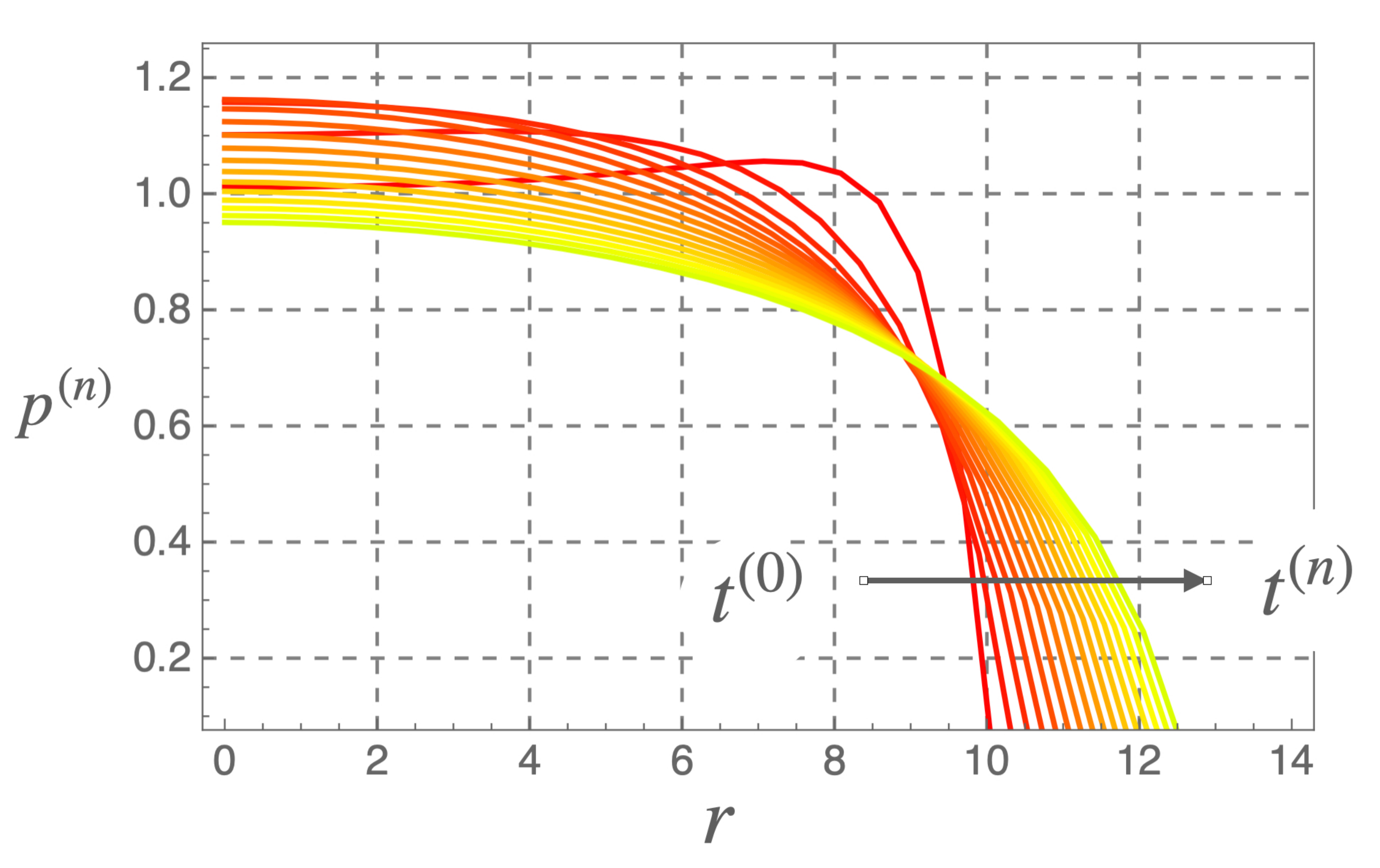}
\caption{}
\end{subfigure}
\begin{subfigure}{0.5\textwidth}
  \includegraphics[width=.9\linewidth]{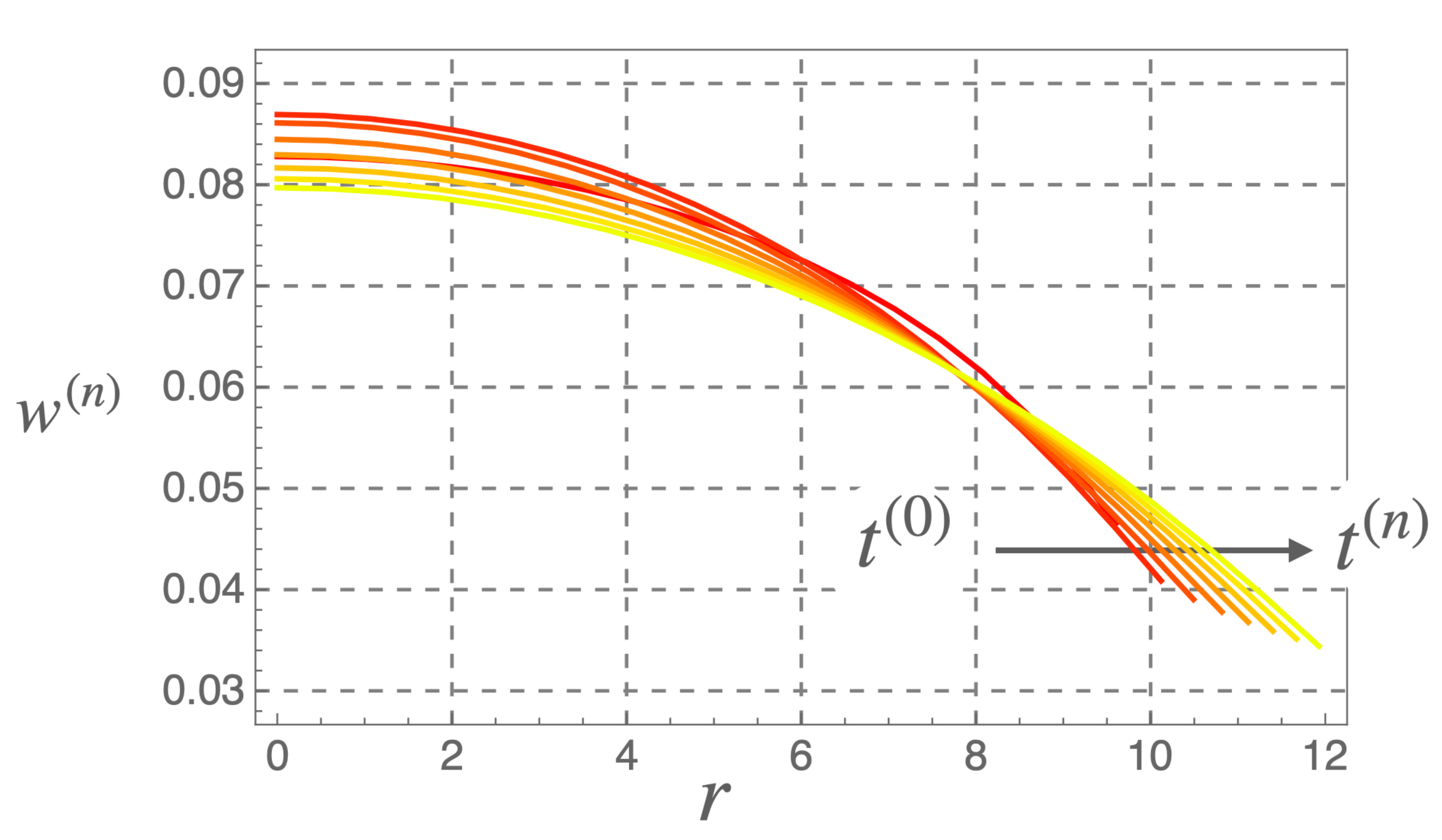}
\caption{}
\end{subfigure}
  \caption{\it Pressure field $p^{(n)}$ (a) and crack opening $w^{(n)}$ (b) in terms of radial coordinate $r$ with $0\leq r \leq \ell(t_n)$ vs time $t_n$. The plane strain modulus of the surrounding medium is equal to $E'=4\cdot 10^2\,MPa.$  }
  \label{Fig::num-sim1}
\end{figure}

\begin{figure}[h]
\begin{subfigure}{0.5\textwidth}
\centering
  \includegraphics[width=.9\linewidth]{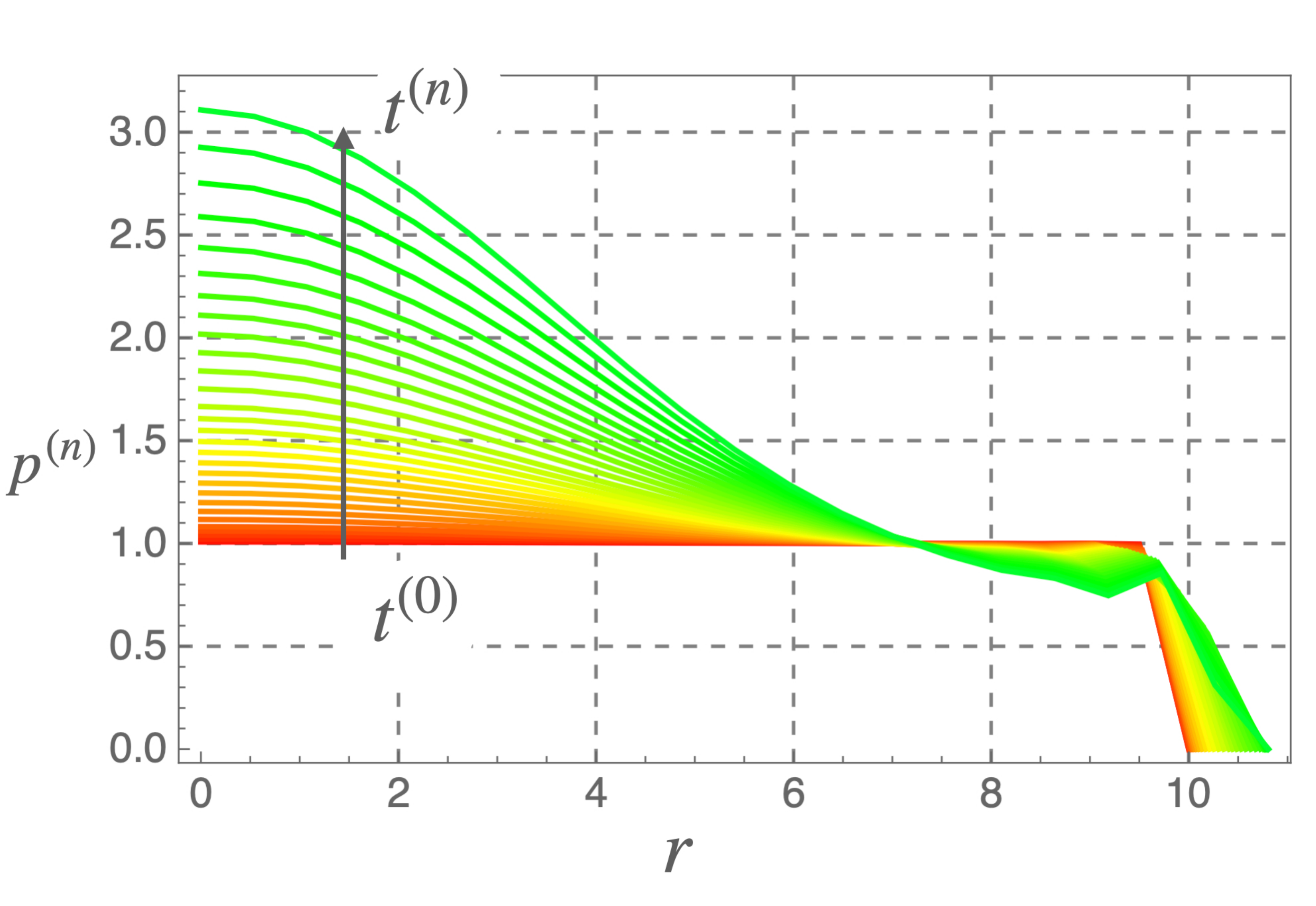}
\caption{}
\end{subfigure}
\begin{subfigure}{0.5\textwidth}
  \includegraphics[width=.9\linewidth]{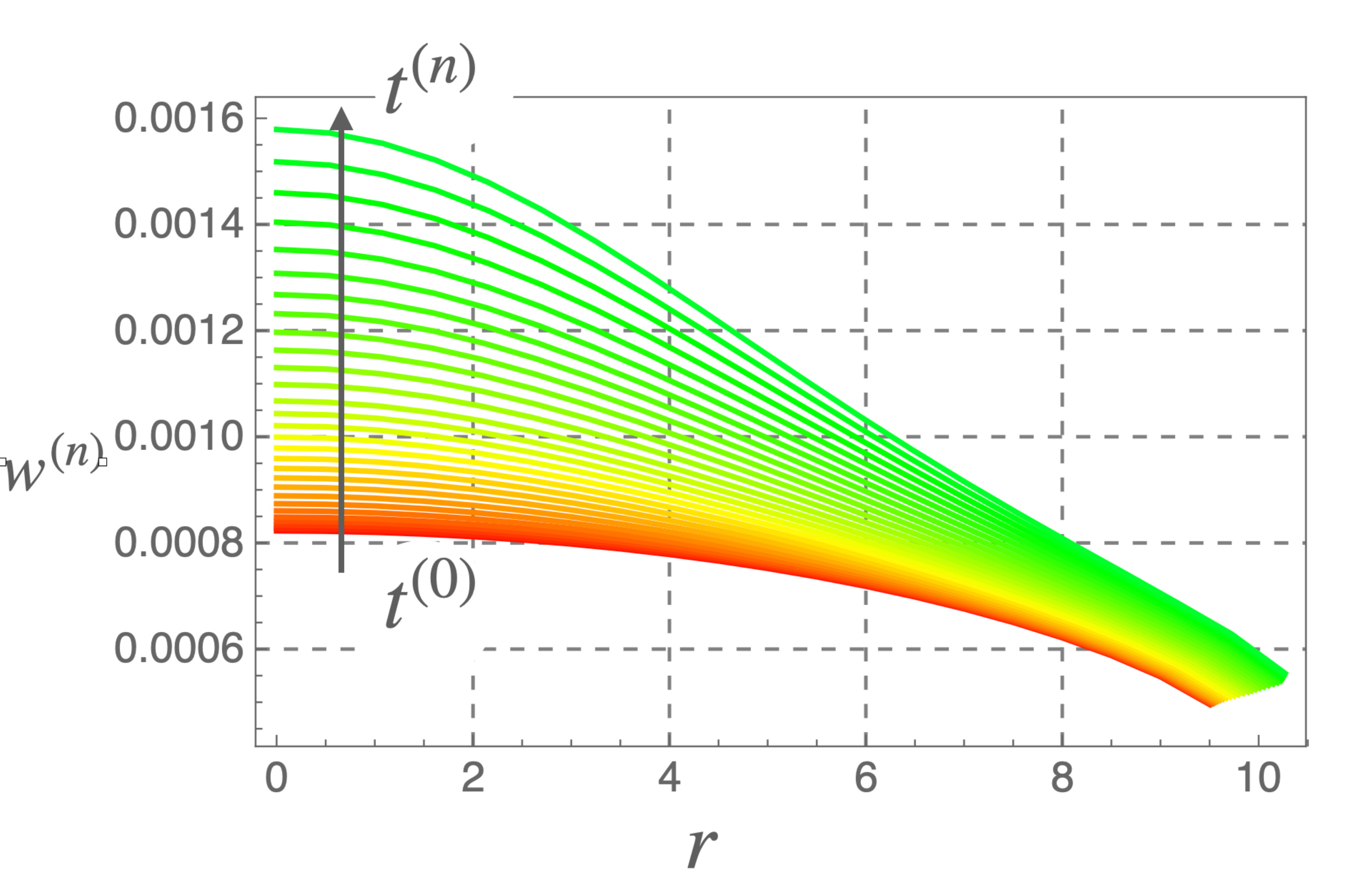}
\end{subfigure}
  \caption{\it Pressure field $p^{(n)}$ (a) and crack opening $w^{(n)}$ (b) in terms of radial coordinate $r$ with $0\leq r \leq \ell(t_n)$ vs time $t_n$. The plane strain modulus of the surrounding medium is equal to $E'=4\cdot10^4\,MPa.$ }
  \label{Fig::num-sim2}
\end{figure}

\begin{figure}[h]
\begin{subfigure}{0.5\textwidth}
\centering
  \includegraphics[width=.9\linewidth]{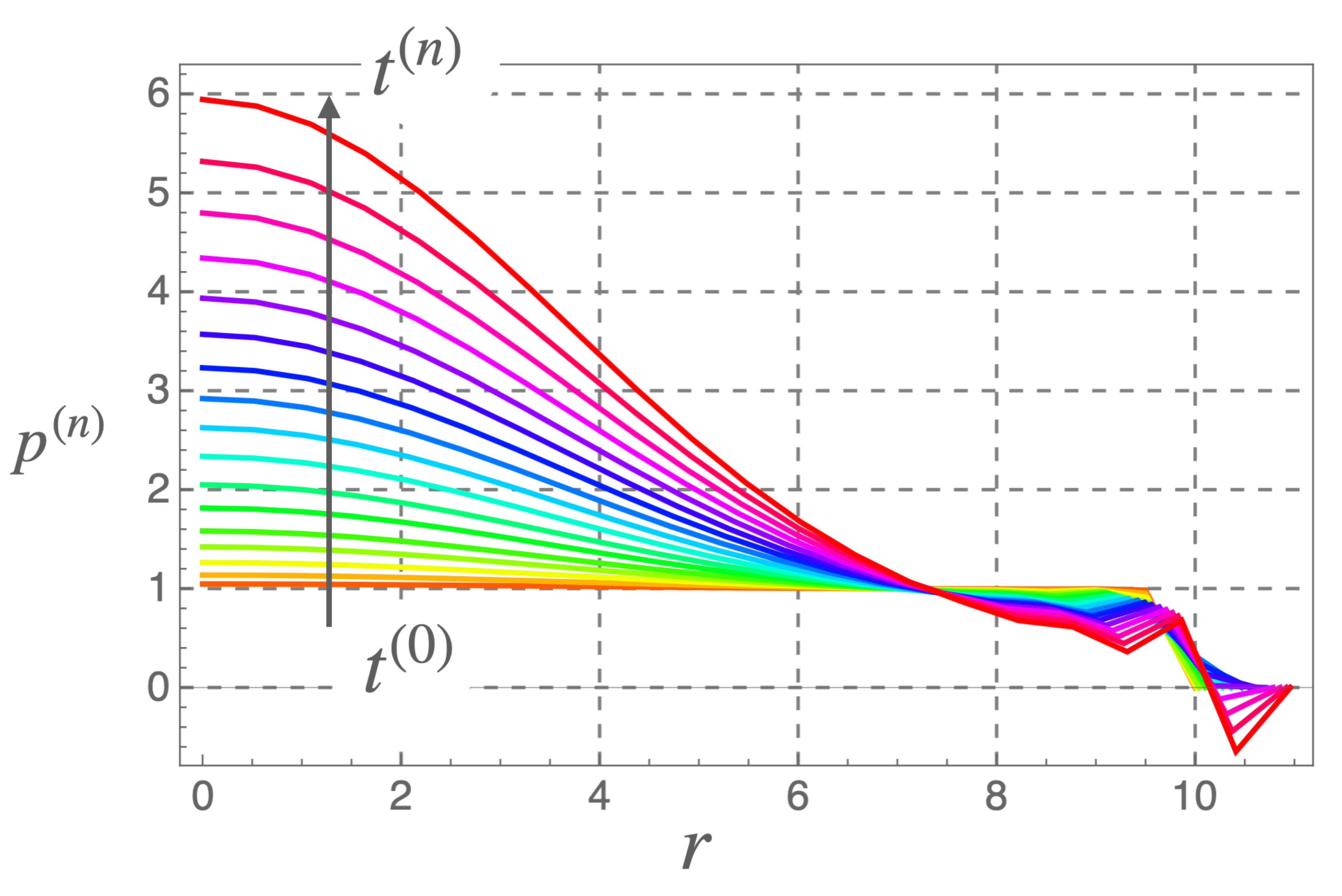}
\caption{}
\end{subfigure}
\begin{subfigure}{0.5\textwidth}
  \includegraphics[width=.9\linewidth]{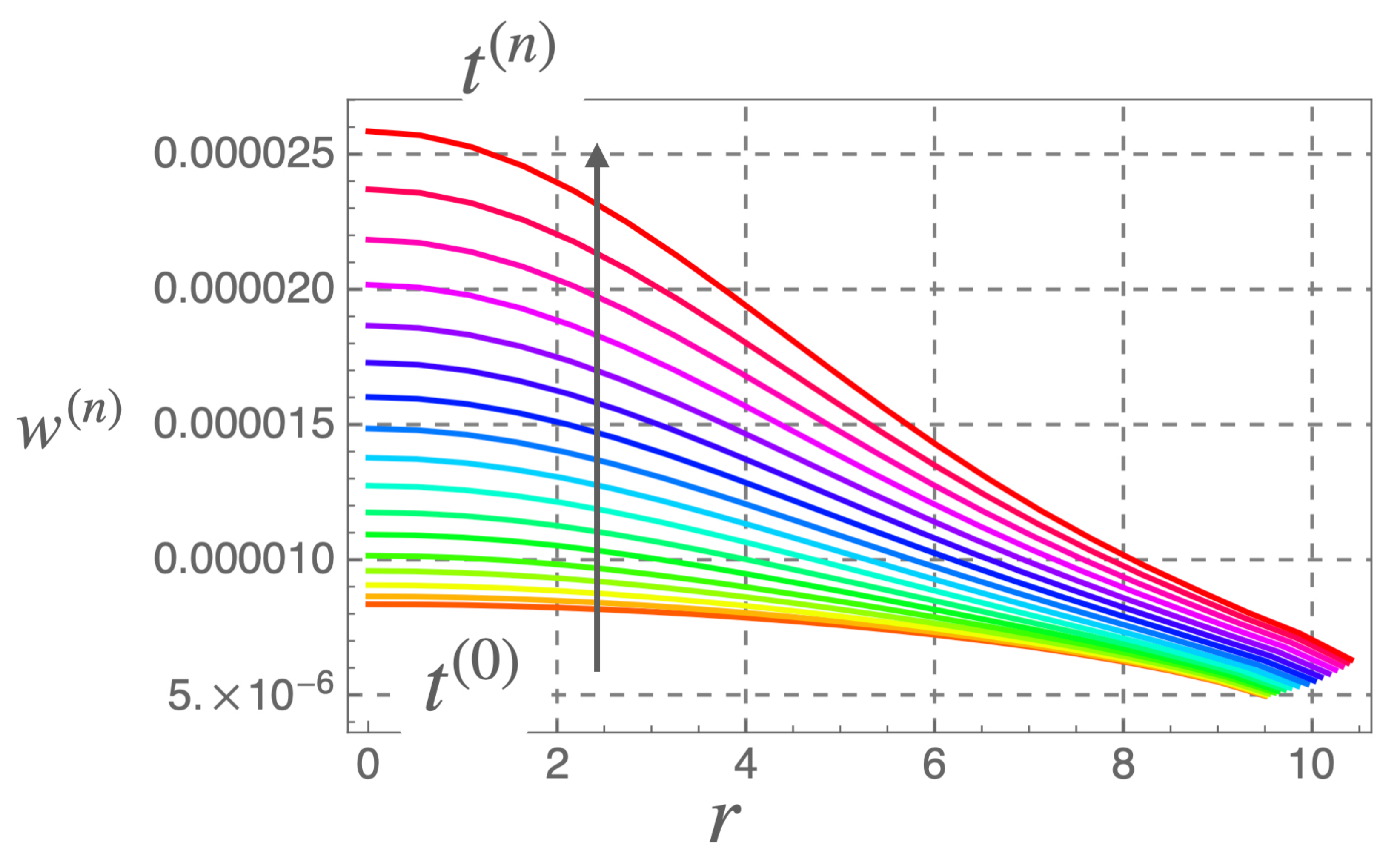}
\caption{}
\end{subfigure}
  \caption{\it Pressure field $p^{(n)}$ (a) and crack opening $w^{(n)}$ (b) in terms of radial coordinate $r$ with $0\leq r \leq \ell(t_n)$ vs time $t_n$. The plane strain modulus of the surrounding medium is equal to $E'=4\cdot10^6\,MPa.$ }
  \label{Fig::num-sim3}
\end{figure}
The net fluid pressure $p^{(n)}(r)$, with $0\leq r\leq \ell(t_n)$, at each time step $t_n = n \, \Delta t$ with a time increment $\Delta t=5\cdot 10^{-7}s$, is shown in Figures \ref{Fig::num-sim1}(a) - \ref{Fig::num-sim3}(a), while Figures \ref{Fig::num-sim1}(b) - \ref{Fig::num-sim3}(b) plot the time evolution of the crack opening $w^{(n)}(r)$ with time $t_n$ in three surrounding media with increasing stiffness.
Pressure increases whereas crack opens less at increasing media stiffness, as expected.
%

At high stiffness, both pressure and opening fields show a strong increment in the correspondence of the inlet, with a steep gradient towards the fluid fronts. The simple numerical analysis here carried out does not allow local mesh refinement, 
which reveals beneficial in the vicinity of $r=\ell^{(n)}$, where oscillations of the solution in the correspondence of the fluid front appear.

In more compliant surrounding medium, the pressure field initially shows a tendency to rise but such a trend is reversed in time, ultimately 
leading to a decrement of both pressure and crack opening, see Fig \ref{Fig::num-sim1}, and to a smoother behavior in space coupled to a more pronounced fluid front advance.
In fact, when $E'\rightarrow 0$, by simple algebra one realizes that given terms of the system (\ref{eq:BackEul}) tend to vanish, and the pressure as well as crack opening along the fluid extent follow the same trend (see equation (\ref{eq:wDisctretization})).
%
%

\section{Conclusions}
\label{Sec::Conclusions}
%

The formulation of a theoretical and computational scheme capable of tracking concurrently two distinct fluid and fracture moving fronts is still a challenging task in hydraulic fracture. Non locality of the equations describing the hydro-mechanical coupling and the lubrication equation describing the fluid response, lead to a non linear system of integro-differential equations for each new geometrical configuration.

The main contribution of the present note stands in the formulation of the crack front tracking, which stems from the fracture propagation depicted as a standard dissipative process. This framework provides analogies with classical viscous regularization approaches in plasticity. They lead to formulate the crack front velocity in an effective form, that enters the differential equations of the hydraulic fracture in a completely novel way. Such an explicit account of the crack front advancing might be a breakthrough in the tracking algorithms currently available in literature.

The set of equations used in this note to build up the algorithm is tailored to an axis-symmetric application used as a proof of concept: the benchmark {shows} that the proposed scheme has great potential and is worth to be extended to general three-dimensional problems, by exploiting either PDEs or integral equations in place of the simple weight function representation formula we made use of in this paper. Discrete weak forms will lead to classical numerical techniques. The requirement of accurate estimation of SIFs along the crack front can be easily achieved by special tip elements as proposed in \cite{ZammarchiEtAlCMAME2017} or in recent literature.

Further extension to permeable surrounding media, where the cavity is filled by a generally unknown pore fluid pressure,  and to different filling media - as gas on non Newtonian fluids -  will be the subject of future research.

\bigskip

\bibliographystyle{unsrt}
\bibliography{Bibliography}


\bigskip \noindent
\textbf{ \large Acknowledgements}

\bigskip
AS, FF gratefully acknowledge the financial support of CeSiA to the project "INTERACTION:INduced ThERmo-chemo-mechanical fRACTure propagatION". AS expresses his  highest gratitude for the appointment of MTS Visiting Professor of Geomechanics at
CEGE, University of Minnesota in year 2019 : this paper is also the outcome of fruitful discussions with E. Detournay and S.Mogilevskaya. Finally, authors are gratefully indebted with M.Jirasek at Czech Technical University, Prague for several enjoyable and fruitful meetings.

\appendix

\section{Appendix A. Lubrication equation proof }
\label{app:DimLubdEq}
The present section is devoted to provide a proof for lubrication equation (\ref{eq:Lubrication1}), starting from Navier-Stokes equations (\ref{eq:NavierStokes}) under lubrication assumptions \cite{Hamrock2004}.
Assuming that the fracture lies in the plane $\{x_1,x_2\}$, for typical lubricant films the ratio between fracture length scale $h$ across its thickness in the $x_3$ direction and lubricant film length scale $L$ in the crack plane  results to be $h/L=O(10^{-3})$.
In order to estimate the orders of magnitude of various terms in the governing equations, dimensionless form of Navier-Stokes    (\ref{eq:NavierStokes}) will be derived in the followings.
 To this aim, one defines  dimensionless coordinates  
 \begin{equation}
 \bar{x}_1=\frac{x_1}{L},
 \hspace{0.2cm}
 \bar{x}_2=\frac{x_2}{L},
 \hspace{0.2cm}
 \bar{x}_3=\frac{x_3}{h},
 \label{eq:AdCoord}
 \end{equation}
which belong to the closed interval $[0,1]$.
Velocity components are adimensionalized as
\begin{equation}
\bar{v}_1=\frac{v_1}{v^*},
\hspace{0.2 cm}
\bar{v}_2=\frac{v_2}{v^*},
\hspace{0.2 cm}
\bar{v}_3=\frac{v_3}{v^*_3},
\label{eq:AdVel}
\end{equation}
with $v^*$ the maximum value of the velocity in the crack plane $\{x_1,x_2\}$, and $v_3^*$ a reference velocity in the $x_3$ direction, which, intuitively, is considerably smaller than $v^*$. 
Furthermore, dimensionless  pressure and mass source are denoted with
\begin{equation}
\bar{p}=Re\,\frac{h}{L}\,\frac{p}{\rho\,{v^*}^2},
\hspace{0.2cm}
\bar{s}=\frac{s}{\rho \,\Omega},
\label{eq:AdQuant}
\end{equation}
where $1/\Omega$ is a characteristic time and $Re=\frac{v^*\,L\,\rho}{\mu}$ is the Reynolds number.
Dimensionless form of   continuity equation (\ref{eq:MassBalancaLocalizedForm2}) thus reads
\begin{equation}
\frac{\partial\bar{v}_1}{\partial \bar{x}_1}+
\frac{\partial\bar{v}_2}{\partial\bar{x}_2}+
\frac{v^*_3}{v^*}\frac{L}{h}\frac{\partial\bar{v}_3}{\partial\bar{x}_3}
=\bar{s}\,\Omega.
\label{eq:AdMassBalance}
\end{equation} 
In order to have all  terms in (\ref{eq:AdMassBalance}) of the same order of magnitude, velocity $v_3^*$ must be equal to $v^* h/L$, meaning that $v^*_3/v^*=O(10^{-3})$.
Plugging relations (\ref{eq:AdCoord}), (\ref{eq:AdVel}), and (\ref{eq:AdQuant}) into (\ref{eq:NavierStokes}), dimensionless form of Navier-Stokes equations results to be
\begin{eqnarray}
\left\{
\begin{array}{l}
-\frac{\partial \bar{p}}{\partial \bar{x}_1}
+
\left(\frac{h}{L}\right)^2\left(\frac{\partial^2\bar{v}_1}{\partial\bar{x}_1^2}+
\frac{\partial^2\bar{v}_1}{\partial\bar{x}_2^2}\right)
+
\frac{\partial^2\bar{v}_1}{\partial\bar{x_3}^2}
+
\frac{\Omega \,h^2}{L\,v^*}\frac{\partial\bar{s}}{\partial\bar{x}_1}=0,
\\
-\frac{\partial \bar{p}}{\partial \bar{x}_2}
+
\left(\frac{h}{L}\right)^2\left(\frac{\partial^2\bar{v}_2}{\partial\bar{x}_1^2}+
\frac{\partial^2\bar{v}_2}{\partial\bar{x}_2^2}\right)
+
\frac{\partial^2\bar{v}_2}{\partial\bar{x_3}^2}
+
\frac{\Omega \,h^2}{L\,v^*}\frac{\partial\bar{s}}{\partial\bar{x}_2}=0,
\\
-\frac{\partial \bar{p}}{\partial \bar{x}_3}+
\left(\frac{h}{L}\right)^4\left(\frac{\partial^2\bar{v}_3}{\partial\bar{x_1}^2}+
\frac{\partial^2\bar{v}_3}{\partial\bar{x_2}^2}\right)+
\left(\frac{h}{L}\right)^2\frac{\partial^2\bar{v}_3}{\partial\bar{x_3}^2}+
\frac{\Omega \,h^2}{L\,v^*}\frac{\partial\bar{s}}{\partial\bar{x}_3}=0.
\end{array}
\right.
\label{eq:AdNavStokes}
\end{eqnarray}
As detailed in section \ref{subsubsec:lubrEq}, the injection of  mass $s$  is uniform across the height of the fracture, meaning that $\partial\bar{s}/\partial\bar{x}_3=0$ in the third of (\ref{eq:AdNavStokes}).
Neglecting terms that are $O(10^{-6})$ and $O(10^{-12})$ in equations (\ref{eq:AdNavStokes}), the only contributions to velocity gradients are $\partial\bar{v}_1/\partial \bar{x}_3$ and $\partial\bar{v}_2/\partial \bar{x}_3$.
Navier-Stokes equation (\ref{eq:NavierStokes}), under lubrication assumptions, eventually simplifies into
\begin{equation}
\left\{
\begin{array}{l}
-\frac{\partial p}{\partial x_1}+\mu\frac{\partial^2 v_1}{\partial x_3^2}+\frac{\mu}{\rho}\frac{\partial s}{\partial x_1}=0,
\\
-\frac{\partial p}{\partial x_2}+\mu\frac{\partial^2 v_2}{\partial x_3^2}+\frac{\mu}{\rho}\frac{\partial s}{\partial x_2}=0,
\\
-\frac{\partial p}{\partial x_3}=0.
\end{array}
\right.
\end{equation}
Double integration of the first and second one leads to
\begin{eqnarray}
v_1&=&-\frac{1}{2\rho}\frac{\partial s}{\partial x_1} x_3^2+
\frac{1}{2\mu}\frac{\partial p}{\partial x_1}x_3^2+
A\,x_3+B,
\nonumber\\
v_2&=&-\frac{1}{2\rho}\frac{\partial s}{\partial x_2} x_3^2+
\frac{1}{2\mu}\frac{\partial p}{\partial x_2}x_3^2+
C\,x_3+D,
\end{eqnarray}
and the value of constants $A, B, C$, and $D$ is computed imposing that $v_1$ and $v_2$ vanish
 at $x_3=-w/2$ and $x_3=w/2$, with $w$ the crack opening.
 Velocity components $v_1$ and $v_2$ eventually read
 \begin{eqnarray}
 v_1&=&\frac{1}{2\mu}\left(x_3^2-\frac{w^2}{4}\right)\left(\frac{\partial p}{\partial x_1}-\frac{\mu}{\rho}\frac{\partial s}{\partial x_1}\right),
 \nonumber\\
  v_2&=&\frac{1}{2\mu}\left(x_3^2-\frac{w^2}{4}\right)\left(\frac{\partial p}{\partial x_2}-\frac{\mu}{\rho}\frac{\partial s}{\partial x_2}\right).
  \label{eq:velocitiescompionents}
 \end{eqnarray}
 Taking into account equations (\ref{eq:velocitiescompionents}), from mass balance equation (\ref{eq:MassBalancaLocalizedForm2}), derivative of component $v_3$ along the thickness can be written as
 \begin{eqnarray}
 \frac{\partial v_3}{\partial x_3}&=&-\frac{\partial v_1}{\partial x_1}-
 \frac{\partial v_2}{\partial x_2}+\frac{s}{\rho}=
 \nonumber\\
 &=&
 \frac{1}{2\mu}
 \left[
 \left(x_3^2-\frac{w^2}{4}
 \right)
 \left(\frac{\mu}{\rho}
 \nabla^2_{\Gamma_w}[s]-
 \nabla^2_{\Gamma_w}[p]
 \right)
 -\frac{w}{2}
 \nabla_{\Gamma_w}[w]\cdot
 \left(
 \frac{\mu}{\rho}
 \nabla_{\Gamma_w} [s]-
 \nabla_{\Gamma_w} [p]
 \right)
 \right]
 +\frac{s}{\rho},
 \label{eq:derivativeofv3}
 \end{eqnarray}
 where
subscript $_{\Gamma_w}$ denotes that the operator to which it refers is  restricted to the fracture plane.
Integration of (\ref{eq:derivativeofv3}) along the thickness finally leads to  expression (\ref{eq:Lubrication1}) of total derivative of fracture opening with respect to time, namely
\begin{equation}
\int_{-w/2}^{w/2}\frac{\partial v_3}{\partial x_3}\,dx_3=\frac{d w}{dt}=\frac{s_w}{\rho_w}
-\surfacedivergence{\frac{w^3}{12\mu_w}
\left(
\frac{\mu_w}{\rho_w}
\nabla_{\Gamma_w}
\left[\frac{s_w}{w}\right ]
-
\nabla_{\Gamma_w}
\left[p_w\right]
\right)
}{\Gamma_w},
\end{equation}
where variables enhanced with subscript $w$ are defined in eq. (\ref{eq:lubr_p_and_s}).

\section{Appendix B. Integral influence functions $A_k(r,t)$}
\label{app:Ak}
One supposes that the one-dimensional domain $[0,\ell(t)]$ is spatially discretized with $Nel$ elements having uniform length $h(t)=\ell(t)/N_{el}$.
If the unknown pressure field  is approximated as $p(r,t)=\sum_{k=0}^{Nel}\psi_k(r)p_k(t)$, where $\psi_k(r)$ are linear shape functions, one has
\begin{equation}
\psi_k(r)=
\left\{
\begin{array}{l}
\psi_k^+(r)=\frac{r}{h}-(k-1) \hspace{0.2cm}if \hspace{0.2cm}(k-1)h\leq r\leq k h \hspace{0.2cm} 
and\hspace{0.2cm} 1\leq k\leq Nel\\
\psi_k^-(r)=-\frac{r}{h}+(k+1) \hspace{0.2cm}if\hspace{0.2cm} kh\leq r\leq (k+1) h \hspace{0.2cm} 
and\hspace{0.2cm} 0\leq k\leq (Nel-1)\\
\end{array}
\right.
\end{equation}
In order to compute discrete integral influence functions $A_k(r,t)$ that relate the crack opening to the pressure field by means of formulas
\begin{equation}
w(r,t)=\sum_{k=0}^{Nel}A_k(r,t)p_k(t),
\end{equation}
and
\begin{equation}
A_k(r,t)=\int_r^{a(t)}
\frac{1}{\sqrt{z^2-r^2}}
\int_0^z
\frac{y \psi_k(y)}{\sqrt{z^2-y^2}}
\,dy\,dz
\label{eq:AkAppendice}
\end{equation}
One has to distinguish among four options, depending on the position of $r$, namely $r<(k-1)h$, $(k-1)h\leq r < k\,h$, $k\, h\leq r < (k+1)h$, or $r\geq (k+1) h$, with $k=0,...,Nel$. Bearing in mind that the elements size is a time dependent quantity, since remeshing has to be performed each time that fluid length $\ell(t)$ varies, time dependence of $A_k$ will be omitted from now on for the sake of readability.
In the first case ($r<(k-1)h$), influence function $A_k$ can be written as
\begin{eqnarray}
&&
A_k(r)=
\int_{(k-1)h}^{k\,h}\frac{1}{\sqrt{z^2-r^2}}I_1(z)\,dz+
\int_{k\,h}^{(k+1)h}\frac{1}{\sqrt{z^2-r^2}}\left(
I_2(z)+I_3(z)
\right)\,dz
+
\nonumber\\
&&+
\int_{(k+1)h}^{a}\frac{1}{\sqrt{z^2-r^2}}\left(
I_2(z)+
I_4(z)
\right)\,dz,
\end{eqnarray}
 where
 \begin{subequations}
 \begin{align}
 I_1(z)&=\int_{(k-1)h}^z
 \frac{y\psi_k^+(y)}{\sqrt{z^2-y^2}}\,dy=\nonumber\\
&= \frac{
 \pi z^2 - 2 h (-1 + k) \sqrt{-h^2 (-1 + k)^2 + z^2} - 
 2 z^2 ArcTan\left(\frac{h (-1 + k)}{\sqrt{-h^2 (-1 + k)^2 + z^2}}\right)}{4 h},\\
I_2(z) &=\int_{(k-1)h}^{k\,h}
 \frac{y\psi_k^+(y)}{\sqrt{z^2-y^2}}\,dy=
 \frac{1}{2 h}
  \left(
  (h-h\,k) \sqrt{-h^2 (-1 + k)^2 + z^2}+
  (- 2h + h\,k) \sqrt{-h^2 k^2 + z^2}+\right.\nonumber\\
  &
    - \left.
   z^2 ArcTan\left(
   \frac{h (-1 + k)}{\sqrt{-h^2 (-1 + k)^2 + z^2}}
   \right) + 
   z^2 ArcTan\left(
   \frac{h k}{\sqrt{-h^2 k^2 + z^2}}\right)
   \right),\\
   I_3(z)&=\int_{k\,h}^{z}
 \frac{y\psi_k^-(y)}{\sqrt{z^2-y^2}}\,dy=\nonumber\\
 &=
 \frac{-\pi z^2 + 2 h (2 + k) \sqrt{-h^2 k^2 + z^2} + 
 2 z^2 ArcTan\left(\frac{h k}{\sqrt{-h^2 k^2 + z^2}}\right)}{4 h},\\
 I_4(z)&=\int_{k\,h}^{(k+1)h}
 \frac{y\psi_k^-(y)}{\sqrt{z^2-y^2}}\,dy=
 \frac{1}{2 h}
  \left(
  (-h-h\,k) \sqrt{-h^2 (1 + k)^2 + z^2}+
  (+ 2h + h\,k) \sqrt{-h^2 k^2 + z^2}+\right.\nonumber\\
  &
    + \left.
   z^2 ArcTan\left(
   \frac{h\,k}{\sqrt{-h^2   k^2 + z^2}}
   \right) - 
   z^2 ArcTan\left(
   \frac{h (k+1)}{\sqrt{-h^2 (k+1)^2 + z^2}}\right)
   \right).
 \end{align}
 \label{eq:internalIntegrals}
 \end{subequations}
When $(k-1)h\leq r < k\,h$ the influence function $A_k(r)$ takes the following form
\begin{eqnarray}
&&
A_k(r)=
\int_{r}^{k\,h}\frac{1}{\sqrt{z^2-r^2}}
\left(I_1(z)-I_1(r)
\right)\,dz 
+
I_1(r)\int_r^{k\,h}\frac{1}{\sqrt{z^2-r^2}}\,dz+
\nonumber\\
&&+
\int_{k\,h}^{(k+1)h}\frac{1}{\sqrt{z^2-r^2}}\left(
I_2(z)+I_3(z)
\right)\,dz
+
\int_{(k+1)h}^{a}\frac{1}{\sqrt{z^2-r^2}}\left(
I_2(z)+
I_4(z)
\right)\,dz,
\label{eq:AKSecondCase}
\end{eqnarray}
where internal integrals $I_1(z),I_2(z),I_3(z)$, and $I_4(z)$ are expressed as in (\ref{eq:internalIntegrals}), and the first two terms on the right hand side of (\ref{eq:AKSecondCase}) take into account that the integrand function of $\int_{r}^{k\,h}\frac{1}{\sqrt{z^2-r^2}}
I_1(z)\,dz $ is singular in $z=r$.
Adopting the same path of reasoning in order to treat singularity points, one has that when $k\,h\leq r <(k+1)k$, $A_k$ can be written as
\begin{eqnarray}
&&
A_k(r)=
\int_{r}^{(k+1)h}\frac{1}{\sqrt{z^2-r^2}}
\left(I_2(z)-I_2(r)+I_3(z)-I_3(r)
\right)\,dz +
\nonumber\\
&&+
\left(I_2(r)+I_3(r)
\right)
\int_r^{(k+1)h}\frac{1}{\sqrt{z^2-r^2}}\,dz
+
\int_{(k+1)h}^{a}\frac{1}{\sqrt{z^2-r^2}}\left(
I_2(z)+I_4(z)
\right)\,dz,
\label{eq:AKThirdCase}
\end{eqnarray}
Finally, when $r\geq (k+1)h$, influence function $A_k(r)$ ca be expressed as
\begin{eqnarray}
&&
A_k(r)=
\int_{r}^{a}\frac{1}{\sqrt{z^2-r^2}}
\left(I_2(z)-I_2(r)+I_4(z)-I_4(r)
\right)\,dz +
\nonumber\\
&&+
\left(I_2(r)+I_4(r)
\right)
\int_r^{a}\frac{1}{\sqrt{z^2-r^2}}\,dz.
\label{eq:AKFourthCase}
\end{eqnarray}
{\itshape Computation of $\frac{\partial A_k(r)}{\partial r}$}
\\
\\
For what regards the computation of the derivative of the influence function $A_k$ with respect to the radial coordinate $r$, when $r< (k-1)h$, being the integration limits independent upon $r$, one simply has
\begin{eqnarray}
&&
\frac{\partial A_k(r)}{\partial r}=
\int_{(k-1)h}^{k\,h}\frac{\partial}{\partial r}
\left(\frac{1}{\sqrt{z^2-r^2}}\right)I_1(z)\,dz+
\int_{k\,h}^{(k+1)h}\frac{\partial}{\partial r}\left(\frac{1}{\sqrt{z^2-r^2}}\right)\left(
I_2(z)+I_3(z)
\right)\,dz
+
\nonumber\\
&&+
\int_{(k+1)h}^{a}\frac{\partial}{\partial r}\left(\frac{1}{\sqrt{z^2-r^2}}\right)\left(
I_2(z)+
I_4(z)
\right)\,dz=
\int_{(k-1)h}^{k\,h}\frac{r}{\sqrt{(z^2-r^2)^3}}
I_1(z)\,dz+\nonumber\\
&&+
\int_{k\,h}^{(k+1)h}\frac{r}{\sqrt{(z^2-r^2)^3}}\left(
I_2(z)+I_3(z)
\right)\,dz
+
\int_{(k+1)h}^{a}\frac{r}{\sqrt{(z^2-r^2)^3}}\left(
I_2(z)+
I_4(z)
\right)\,dz.
\end{eqnarray}
In the second case, when $(k-1)h\leq r < k\,h$, the expression of $\frac{\partial A_k(r)}{\partial r}$ becomes
\begin{eqnarray}
&&
\frac{\partial A_k(r)}{\partial r} = \lim_{z\rightarrow r^+}\frac{I_1(z)-I_1(r)}{\sqrt{z^2-r^2}}
+
\int_{r}^{k\,h}
\frac{\partial}{\partial r}
\left(
\frac{I_1(z)-I_1(r)}{\sqrt{z^2-r^2}}
\right)\,dz
+
\frac{\partial I_1(r)}{\partial r}
\int_r^{k\,h}
\frac{1}{\sqrt{z^2-r^2}}
+
\nonumber\\
&&
+
I_1(r)
\frac{\partial}{\partial r}
\int_r^{k\,h}
\frac{1}{\sqrt{z^2-r^2}}\,dz
+
\int_{k\,h}^{(k+1)h}
\frac{r}{\sqrt{(z^2-r^2)^3}}(I_2(z)+I_3(z))\,dz
+
\nonumber\\
&&
+
\int_{(k+1)h}^a
\frac{r}{\sqrt{(z^2-r^2)^3}}(I_2(z)+I_4(z))\,dz,
\label{eq:dAkdrSecondCase}
\end{eqnarray}
where the limit in the right hand side of (\ref{eq:dAkdrSecondCase}) vanishes, and
\begin{eqnarray}
\frac{\partial I_1(r)}{\partial r}
=
\frac{
r \left(
\pi - 2 
ArcTan
\left(
\frac{ h (-1 + k)}{\sqrt{-h^2 (-1 + k)^2 + r^2}}
\right)\right)}{2 h},
\hspace{0.5cm}
\frac{\partial}{\partial r}\int_r^{k\,h}
\frac{1}{\sqrt{z^2-r^2}}\,dz=
-\frac{ k\,h}{r \sqrt{h^2 k^2 - r^2}}
\end{eqnarray}
In the third case, when r $k\,h\leq r <(k+1) h$, the derivative of the influence function $A_k$ with respect to $r$ can be written as
\begin{eqnarray}
&&\frac{\partial A_k(r)}{\partial r}=
\lim_{z\rightarrow r^+}\frac{I_2(z)+I_3(z)-I_2(r)+I_3(r)}{\sqrt{z^2-r^2}}
+
\int_r^{(k+1)h}
\frac{\partial}{\partial r}
\left(
\frac{I_2(z)+I_3(z)-I_2(r)+I_3(r)}{\sqrt{z^2-r^2}}
\right)\,dz
+
\nonumber\\
&&+
\frac{\partial (I_2(r)+I_3(r))}{\partial r}
\int_r^{(k+1) h}\frac{1}{\sqrt{z^2-r^2}}\,dz
+
(I_2(r)+I_3(r))
\frac{\partial}{\partial r}
\left(
\int_r^{(k+1) h}\frac{1}{\sqrt{z^2-r^2}}\,dz
\right)+
\nonumber\\
&&+
\int_{(k+1)h}^a
\frac{r}{\sqrt{(z^2-r^2)^3}}(I_2(z)+I_4(z))\,dz,
\label{eq:dAkdrThirdCase}
\end{eqnarray}
where, again, the limit on the right hand side of (\ref{eq:dAkdrThirdCase}) vanishes, and
\begin{eqnarray}
&&
\frac{\partial I_2(r)}{\partial r}
=
r \left(
-
\frac{1}{\sqrt{-h^2 k^2 + r^2}}
-
\frac{ ArcTan\left(\frac{h (-1 + k)}{\sqrt{-h^2 (-1 + k)^2 + r^2}}\right)}{h}
+
\frac{ ArcTan\left(\frac{h k}{\sqrt{-h^2 k^2 + r^2}}\right)}{h}
\right),
\nonumber\\
&&
\frac{\partial I_3}{\partial r}
=
-\frac{\pi r}{2 h}
+\frac{r}{\sqrt{-h^2 k^2 + r^2}} 
+\frac{r ArcTan\left(\frac{h k}{\sqrt{-h^2 k^2 + r^2}}\right)}{h},
\nonumber\\
&&
\frac{\partial}{\partial r}
\left(
\int_r^{(k+1) h}\frac{1}{\sqrt{z^2-r^2}}\,dz
\right)=-\frac{(k+1)h}{r \sqrt{h^2 (k+1)^2 - r^2}}.
\end{eqnarray}
Finally, when $r\geq (k+1) h$, one has
\begin{eqnarray}
&&
\frac{\partial A_k(r)}{\partial r}=
\lim_{z\rightarrow r^+}
\frac{I_2(z)+I_4(z)-I_2(r)-I_4(r)}{\sqrt{z^2-r^2}}
+
\int_r^a
\frac{\partial}{\partial r}
\left(
\frac{I_2(z)+I_4(z)-I_2(r)-I_4(r)}{\sqrt{z^2-r^2}}
\right)\,dz+
\nonumber\\
&&
+
\frac{\partial(I_2(r)+I_4(r))}{\partial r}
\int_r^a
\frac{1}{\sqrt{z^2-r^2}}\,dz
+
(I_2(r)+I_4(r))
\frac{\partial}{\partial r}
\left(
\int_r^a
\frac{1}{\sqrt{z^2-r^2}}\,dz
\right),
\label{eq:dAkdrFourthCase}
\end{eqnarray}  
with a vanishing value of the limit in the right hand side of (\ref{eq:dAkdrFourthCase}) and
\begin{eqnarray}
&&
\frac{\partial I_4(r)}{\partial r}=
r 
\left(
\frac{1}{\sqrt{-h^2 k^2 + r^2}}
+\frac{ ArcTan\left(\frac{h k}{\sqrt{-h^2 k^2 + r^2}}\right)}{h} 
- \frac{ArcTan\left(\frac{(k+1)h}{\sqrt{-h^2 (k+1)^2 + r^2}}\right)}{h}
\right),
\nonumber\\
&&
\frac{\partial}{\partial r}
\left(
\int_r^a
\frac{1}{\sqrt{z^2-r^2}}\,dz
\right)=-\frac{a}{r \sqrt{a^2 - r^2}}
\end{eqnarray}
{\itshape Computation of $\frac{\partial A_k(r)}{\partial a}$}
\\
\\
Based upon the definition (\ref{eq:AkAppendice}) of influence function $A_k(r)$, the computation of its derivative with respect to crack radius is quite simple. In fact, being the integral kernel  independent upon $a$, when integration limits are different from $a$, $\frac{\partial A_k(r)}{\partial a}=0$, and when the superior integration limit is equal to $a$, the derivative of $A_k(r)$ with respect to $a$ is equal to the integral kernel evaluated at $z=a$.
Therefore, for all the four cases based upon the position of $r$, one has the same expression for  $\frac{\partial A_k(r)}{\partial a}$ that is equal to
\begin{eqnarray}
&&\frac{\partial A_k(r)}{\partial a}=
\frac{I_2(a)+I_4(a)}{\sqrt{a^2-r^2}}=
-\frac{1}{2 h \sqrt{a^2-r^2}}
\left(
h (k+1) \sqrt{-h^2 (k-1)^2 + a^2}
-
2 h k \sqrt{-h^2 k^2 + a^2}
+\right.\nonumber\\
&&\left.+
h (k+1) \sqrt{-h^2 (k+1)^2 + a^2}
+
a^2 ArcTan\left(\frac{h (k-1)}{\sqrt{-h^2 (k-1)^2 + a^2}}\right)
-
2 a^2 ArcTan\left(\frac{h k}{\sqrt{-h^2 k^2 + a^2}}\right)
+
\right.
\nonumber\\
&&\left.+
a^2 ArcTan\left(\frac{h (k+1)}{\sqrt{-h^2 (k+1)^2 + a^2}}\right)
\right)
\end{eqnarray}
{\itshape Computation of $\frac{\partial A_k(r)}{\partial \ell}$}
\\
\\
The derivative of influence function $A_k(r)$ with respect to the fluid extent $\ell$ is non vanishing only when $k=Nel$.
It has the following expression, both when $r< (N_{el}-1)h$ and when $(N_{el}-1)h\leq r< N_{el}h$
\begin{eqnarray}
\frac{\partial A_{N_{el}}(r)}{\partial \ell}=\frac{I_1(\ell)}{\sqrt{\ell^2-r^2}}=
\frac{-2 ( N_{el}-1) \sqrt{h^2 (-1 + 2 N_{el})}
+h N_{el}^2 \pi
+
-2 h N_{el}^2 ArcTan\left(
\frac{h (N_{el}-1}{\sqrt{h^2 (-1 + 2 N_{el}}}\right)
}{4\sqrt{\ell^2-r^2}}
\end{eqnarray}

\end{document}